# Improved adjoint lattice Boltzmann method for topology optimization of laminar convective heat transfer

Ji-Wang Luo [a, b], Li Chen [a, *], Kentaro Yaji [b, *], Wen-Quan Tao [a]

a: Key Laboratory of Thermo-Fluid Science and Engineering of MOE, School of Energy and Power Engineering, Xi'an Jiaotong University, Xi'an, Shaanxi 710049, China b: Department of Mechanical Engineering, Graduate School of Engineering, The University of Osaka, 2-1, Yamadaoka, Suita, Osaka 565-0871, Japan

*Corresponding author:

Li Chen Email: lichennht08@mail.xjtu.edu.cn

Kentaro Yaji Email: yaji@mech.eng.osaka-u.ac.jp

**Abstract**

Solving flow-related inverse problems such as topology optimization problems is intricate but significant in various engineering fields. The lattice Boltzmann method (LBM) and the related adjoint method are highly suitable to perform sensitivity analysis in flow-related inverse problems thanks to their strong capability to handle complex structures and excellent parallel scalability. However, the current continuous adjoint LBM shows theoretical inconsistency and poor numerical stability for open flow systems. To solve these issues, the present work develops the fully consistent adjoint boundary conditions from the discrete adjoint LBM. For the first time, the gap between the two adjoint LBMs is unveiled by rigorously deriving both the continuous and discrete adjoint LBMs and comprehensively evaluating their numerical performances in the 2D and 3D pipe bend optimization cases. It is revealed that theoretical inconsistency or singularity exists in the continuous adjoint boundary conditions for open flow systems, corresponding to a much inferior numerical stability of the adjoint solution and an obvious numerical error in sensitivity. Fully consistent adjoint boundary conditions in elegant local form are derived from the discrete adjoint LBM in this work, which can always acquire exact sensitivity results and the theoretically highest numerical stability, with a 10 times higher Reynolds number (*Re*) achieved while without any increase of computational cost. 3D microchannel heat sinks under various *Re* are designed, and the esthetic and physically reasonable optimized designs are obtained under various parameter settings, demonstrating the necessity and versatility of the presented discrete adjoint LBM.

## 1 Introduction

Solving the flow-related inverse problems, such as flow identification problem [1], structural optimization problem [2] or flow control problem [3], is intricate but significant in various engineering fields. One of the recent research hotspots is the flow-related topology optimization (TO). Although TO was first developed and matured in the field of solid mechanics [4], its applications to fluid flow problem have been proceeded rapidly, from the creeping Stokes flow [5], the steady-state [6] and unsteady-state [7] Navier-Stokes flow to turbulent flow [8] and rarefied flow [9]. Besides, the achievements of TO in the flow-related multi-physics processes are also noteworthy, for example, TO has been applied to design various kinds of heat sinks [2], heat exchangers [10], chemical reactors [11] and porous electrodes [12]. In the flow-related inverse problems like TO, usually there are large numbers of design or control variables to be optimized. For such large-scale optimization problem, usually the gradient-based algorithm is adopted to find an optimal solution, and it requires an additional sensitivity analysis to obtain the gradient or sensitivity [13]. To develop robust, accurate and efficient numerical methods for the sensitivity analysis is one of the major research tasks of the research community.

The sensitivity can be obtained by various methods, such as the finite difference method (FDM), complex differentiation method and adjoint method, among which the adjoint method is the most popular due to its low computational cost [13]. The adjoint method introduces the auxiliary adjoint variables to simplify the calculation of sensitivity. By solving the adjoint equations, whose computational cost is almost equal to that of one physical simulation (also called the primal or forward simulation), the sensitivity can be readily obtained. Thanks to its high computational efficiency, the adjoint method has been widely adopted to perform the sensitivity analysis in solving the flow-related inverse problems, especially in the field of flow-related TO [14]. Therefore, to perform efficient sensitivity analysis, both the forward and adjoint solvers need to be robust, accurate and efficient. To simulate the flow-related processes, usually the macroscopic discretization methods such as the finite volume method (FVM) and the finite element method (FEM) are adopted. Recently, the lattice Boltzmann method (LBM) based on the kinetic theory has become a promising alternative outside the FVM and FEM. The LBM has been applied to solve various intricate flow problems such as multi-phase flow and porous flow, thanks to its strong capability to deal with complicated interfaces and its extraordinary parallel scalability [15-17]. Besides, the LBM has also been applied to solve the flow-related inverse problems, especially the flow-related TO problems, by developing the corresponding adjoint LBM [18].

Generally, the adjoint method can be classified into the discrete or continuous adjoint method, in which the discrete adjoint method follows the idea of 'first discretize then adjoint' while the continuous one takes the way of 'first adjoint then discretize' [13]. In the field of adjoint LBM, the discrete adjoint LBM was first developed by Tekitek et al. [19] to solve a parameter identification problem. Then Pingen et al. [20] further developed the discrete adjoint LBM for solving the flow-related TO problems, such as the TO problems involving steady-state Navier-Stokes flow [20] and non-Newtonian flow [21]. Kreissl et al. [22] extended Pingen's approach to the TO study about the fluid-structure interaction problems, and they also integrated the discrete adjoint LBM with explicit level set approach [23], while Makhija et al. [24] applied the similar method to design the micro-mixers considering convective mass transfer. It is worth pointing out that in the above early studies adopting the discrete adjoint LBM (except Tekitek et al. [19]), the Jacobian of discretized system, which is usually a large-scale matrix, has to be composed and inverted, which goes against the simplicity and locality of the primal LBM and makes the discrete adjoint LBM almost infeasible for large-scale optimization problems. Later, Krause and Heuveline [25] presented the combination of discrete adjoint LBM with automatic differentiation techniques. Liu et al. [26] presented the discrete adjoint LBM with the multiple-relaxation time (MRT) collision operator, and they circumvented the use of Jacobian by separately handling the collision and streaming operators, similar to the idea of Tekitek [19]. Such strategy brings the adjoint solution process back to the way of simplicity and locality, greatly improving the solution efficiency. Therefore, almost all of the subsequent studies on the discrete adjoint LBM adopted that strategy. Hekmat and Mirzaei [27] developed the discrete adjoint LBM based on both macroscopic and microscopic adjoint concepts, and they also compared the discrete and continuous adjoint LBM regarding the numerical accuracy and convergence behavior [28]. Vergnault and Sagaut [29] developed the discrete adjoint LBM for noise control problems. Yonekura and Kanno [30] proposed the method to simplify the calculation of sensitivity by only using the transient information. Łaniewski-Wołłk and Rokicki [31] developed the discrete adjoint LBM by combining hand derivation and automatic differentiation to solve TO problems of conjugate heat transfer. Nørgaard et al. [32] developed the unsteady-state discrete adjoint LBM for flow control and pump design problems, as well for pressure diode problem [33]. Recently, Cheylan et al. [34] presented the shape optimization of car in turbulent flow using the discrete adjoint LBM with two-relaxation-time collision model. Rutkowski et al. [3] developed the discrete adjoint LBM for open-loop control of flapping wing motion. Zarth et al. [35] presented the shape optimization for fluid flow using discrete adjoint LBM combined with automatic differentiation. Kusano [36] proposed the sensitivity analysis for flow-induced sound problems based on the discrete adjoint LBM.

Along with the progresses in the discrete adjoint LBM, the continuous adjoint LBM also advanced rapidly. By developing the continuous adjoint LBM, deeper mathematical and physical understanding to the adjoint problems can be obtained, and the more flexible and cheaper solution strategies can be developed [13]. In 2013, Krause et al. [37] presented the flow control and optimization using the continuous adjoint LBM, which was later extended by Yaji et al. [38] to the TO of fluid flow problems. Yaji et al. [39] further developed the continuous adjoint LBM based on the discrete velocity Boltzmann equation (DVBE), aiming to introducing the commonly-used boundary schemes into the adjoint derivation. Here the adjoint method based on DVBE is categorized as a continuous adjoint method considering that both the DVBE and its adjoint are continuous in time and space, which need further discretization before the solution procedure, and the discreteness of the DVBE in the velocity space does not alter its essence of 'first adjoint then discretize'. Chen et al. [40] and Nguyen et al. [41] extended such approach to unsteady-state fluid flow problems, and Yaji et al. [42] also applied it to unsteady-state convective heat transfer cases. Dugast et al. [11, 43] developed the continuous adjoint LBM combined with the level set method for convective heat and mass transfer problems. Klemens et al. [1] presented the applications of continuous adjoint LBM to flow characterization and domain identification problems, as well as to the noise reduction for MRI measurements [44]. Li et al. [45] presented the continuous adjoint LBM for airfoil design optimization problems. Khouzani and Moghadam [46] proposed the continuous adjoint LBM based on the circular function scheme. Xie et al. [47] developed the continuous adjoint lattice kinetic scheme for flow optimization, and Tanabe et al. [48] recently extended this approach to thermal flow problems. Luo et al. [2, 14] developed the continuous adjoint LBM to design solid and porous heat sinks in natural convection system, and Tanabe et al. [49] further considered the transient situation. Yong and Zhao [50] developed a low-storage continuous adjoint LBM for transient flow control problems. Khouzani and Moghadam [51] developed the continuous adjoint LBM for laminar compressible flow, and applied to airfoil shape optimization.

It can be concluded from the above literature review that both the discrete and continuous adjoint LBM have been widely concerned and developed in recent years, and their potentials to solve various flow-related inverse problems have been explored. Despite that, one of the key components in the adjoint LBM, the adjoint boundary condition, has rarely been the research emphasis. For example, Hekmat and Mirzaei [27, 28] developed the discrete adjoint LBM and performed the comparison between discrete and continuous adjoint LBM with the simple bounce-back and periodic boundaries considered only, while Liu et al. [26] used the posteriori estimated boundary schemes in the discrete adjoint LBM. Recent works presented various approaches to reach the discrete adjoint LBM, yet still there lacks clear and accurate descriptions to the related adjoint boundary conditions. In the continuous adjoint, the issues pertaining to the adjoint boundary condition are even more severe. When the fully continuous Boltzmann-BGK (Bhatnagar-Gross-Krook) equation is used, the well-established boundary schemes in LBM, like the Zou-He scheme [52], cannot be incorporated into the adjoint derivation, leading to poorer numerical performance in either the primal or the adjoint LBM solution. The continuous adjoint LBM derived from the DVBE allows the use of Zou-He scheme, while the inconsistency still exists in the corresponding adjoint boundary condition [39]. To solve this issue, some assumptions or simplifications to the adjoint boundary conditions have been proposed [45, 48], which renders the adjoint solution possible but comes at the cost of losing numerical accuracy in sensitivity. Besides, poor numerical stability of the resulting continuous adjoint LBM has also been reported [43], that is, the adjoint solution gets diverged once the Reynolds number (*Re*) approaches nearly 25, despite that the primal LBM simulation can give converged solution at a much higher *Re* (e.g., 200), while the causes and solutions of such poor stability remain unclear. Therefore, specialized research to the adjoint boundary conditions in the continuous and discrete adjoint LBM is highly desired.

This work aims to develop robust, accurate and efficient adjoint LBM for the sensitivity analysis of flow-related processes, with a special focus on the adjoint boundary conditions. Specifically, both the continuous and discrete adjoint LBMs will be rigorously derived, with their commonalities and discrepancies revealed and with the clear and accurate descriptions to the resulting adjoint boundary conditions presented. The numerical performances of the two adjoint LBMs, including the numerical accuracy, numerical stability and computational efficiency, are comprehensively evaluated. The inconsistency in the continuous adjoint boundary conditions and its contributions to the numerical instability of the adjoint LBM are well revealed and demonstrated for the first time, and solution to these issues are given from the view of discrete adjoint LBM. Besides, the application of the developed adjoint LBM to the TO of forced convection microchannel heat sink is also presented.

The remainders of the present work are arranged as follows. Section 2 introduces the LBM, the TO method, and derivations of continuous and discrete adjoint LBM. Section 3 performs the comparative study about the numerical performances of the two adjoint LBM, and the TO of 3D microchannel heat sink is also presented. Finally the important conclusions of the present work are drawn in Section 4.

## 2 Methodology

This section will present the derivations of both continuous and discrete adjoint LBM, with the inconsistency issue in the continuous adjoint boundary condition and the differences between the two adjoint models being revealed. Here the considered optimized problem will be the typical TO problem, namely the pipe bend problem, as schematically shown in Fig. 1.

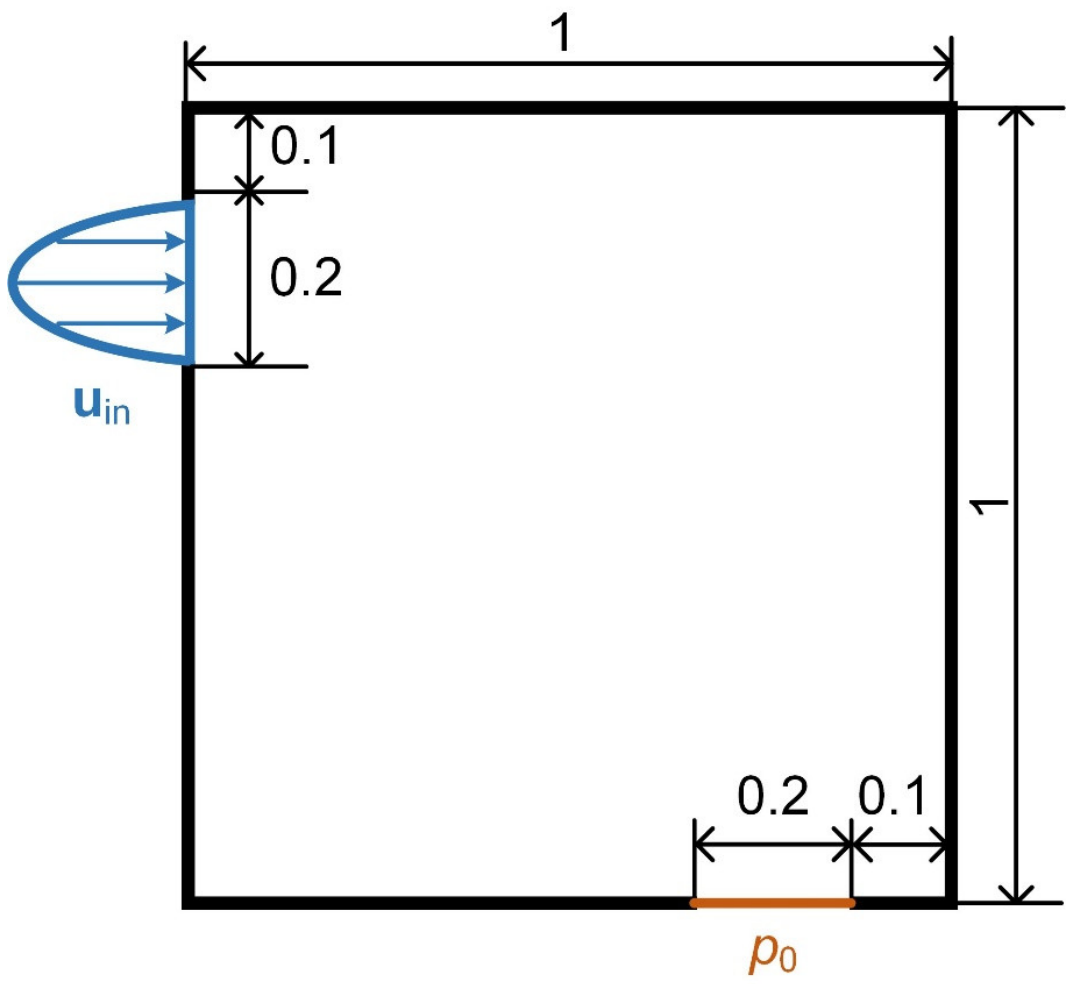

**Figure 1.** Schematic of the 2D pipe bend case.

*2.1 Lattice Boltzmann model*

Since the LBM is used to resolve the fluid flow, the governing equation in the continuous form will be the following discrete velocity Boltzmann equation (DVBE) [53]

$$\frac{\partial f_i}{\partial t} + \mathbf{e}_i \cdot \nabla f_i = -\frac{1}{\tau_c}\left(f_i - f_i^{\text{eq}}\right) \tag{1}$$

where $f_i$ is the distribution function with velocity $\mathbf{e}_i$, at position $\mathbf{x}$ and time $t$. $\tau_c$ is the relaxation time. To numerically solve the above partial differential equation (PDE), one should further discretize it in time and space, leading to the lattice Boltzmann equation (LBE) [53]

$$f_i(\mathbf{x} + \mathbf{e}_i \Delta t, t + \Delta t) - f_i(\mathbf{x}, t) = -\frac{1}{\tau}\left[f_i(\mathbf{x}, t) - f_i^{\text{eq}}(\mathbf{x}, t)\right] \tag{2}$$

where $\Delta t$ is the time step and $\mathbf{e}_i$ is the set of discrete velocities, and $\tau=\tau_c/\Delta t$ is the dimensionless relaxation time. For fluid flow simulation, the D2Q9 stencil and the D3Q19 stencil are the most commonly used discrete velocity sets, as shown in Fig. 2. The dimensionless relaxation time $\tau$ is dependent on the kinetic viscosity $\nu$

$$\nu = \frac{c^2}{3}(\tau - 0.5)\Delta t \tag{3}$$

where $c=\Delta x/\Delta t=1$ in the present work. $f_i^{\text{eq}}$ in Eqs. (1) and (2) is the equilibrium distribution function expressed by [18]

$$f_i^{\text{eq}} = \omega_i\left\{\rho + \rho_0\left[3(\mathbf{e}_i \cdot \mathbf{u}) + \frac{9}{2}(\mathbf{e}_i \cdot \mathbf{u})^2 - \frac{3}{2}\mathbf{u}^2\right]\right\} \tag{4}$$

where $\omega_i$ is the weighting coefficient. Here the incompressible LB model instead of the standard

LB model is used to eliminate the compressibility effect [43]. For D2Q9 and D3Q19 stencils, the values of $\omega_i$ will be

$$\text{D2Q9:}\ \omega_i=\begin{cases}\frac{4}{9}, & i=0\\ \frac{1}{9}, & i=1-4\\ \frac{1}{36}, & i=5-8\end{cases}\qquad \text{D3Q19:}\ \ \omega_i=\begin{cases}\frac{1}{3}, & i=0\\ \frac{1}{18}, & i=1-6\\ \frac{1}{36}, & i=7-18\end{cases}\tag{5}$$

where $\rho_0$ is the constant density. $\rho$ and $\mathbf{u}$ are respectively the density and velocity, which can be calculated by the zeroth and first order moments of $f_i$

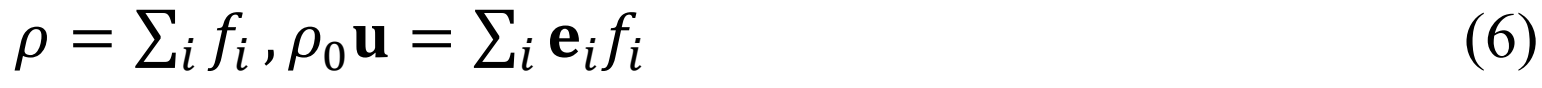

$$\rho=\sum_i f_i\,,\rho_0\mathbf{u}=\sum_i \mathbf{e}_i f_i\tag{6}$$

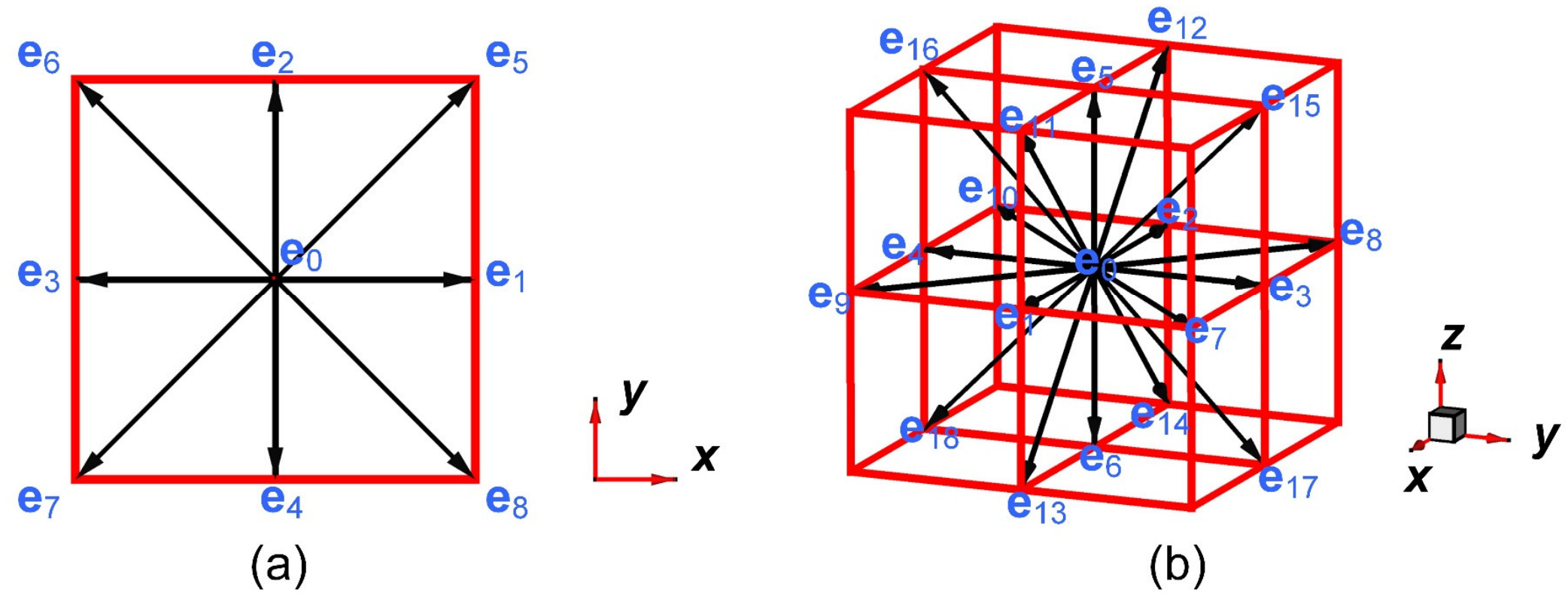

**Figure 2.** Discrete velocities of (a) D2Q9 stencil and (b) D3Q19 stencil.

The unknown distribution functions at the boundary should be calculated by some specific boundary schemes. In this work, the Zou-He boundary scheme, or namely the non-equilibrium bounce-back scheme, is applied for open boundaries [52]. Taking the 2D case as example, the unknown distributions at the left inlet with given velocity $\mathbf{u}_{\text{in}}=(u_{\text{in}}, 0)$ are given by [52]

$$\begin{aligned}
&\rho=f_0+f_2+f_4+2(f_3+f_6+f_7)+\rho_0 u_{\text{in}}\\
&f_1=f_3+\tfrac{2}{3}\rho_0 u_{\text{in}}\\
&f_5=f_7+\tfrac{1}{6}\rho_0 u_{\text{in}}+\tfrac{1}{2}(f_4-f_2)\\
&f_8=f_6+\tfrac{1}{6}\rho_0 u_{\text{in}}-\tfrac{1}{2}(f_4-f_2)
\end{aligned}\tag{7}$$

While those unknowns at the bottom outlet with fixed density (i.e., fixed pressure) $\rho=\rho_0$ can be calculated by [52]

$$\begin{aligned}
&v=1-\big(f_0+f_1+f_3+2(f_4+f_7+f_8)\big)/\rho_0\\
&f_2=f_4+\tfrac{2}{3}\rho_0 v\\
&f_5=f_7+\tfrac{1}{6}\rho_0 v+\tfrac{1}{2}(f_3-f_1)\\
&f_6=f_8+\tfrac{1}{6}\rho_0 v-\tfrac{1}{2}(f_3-f_1)
\end{aligned}\tag{8}$$

Other boundaries are no-slip walls, which is easily realized by the bounce-back scheme

$$f_{i,(\mathbf{e}_i\cdot\mathbf{n}<0)} = f_{\bar{i},(\mathbf{e}_{\bar{i}}\cdot\mathbf{n}>0)} \tag{9}$$

In summary, by sequentially performing the collision, streaming and boundary condition steps, the fluid flow inside the open chamber shown in Fig. 1 can be resolved.

*2.2 Level set-based topology optimization*

The purpose of the pipe bend problem is to design the pipe in an open chamber (Fig. 1) to minimize the dissipation power subject to a constraint of fluid volume fraction. Therefore, the optimization problem could be formulated as follows

$$\begin{aligned} \min_{\alpha} \quad & J = p_{\text{in}} = \frac{1}{3W_{\text{in}}t_{\text{f}}}\int_t \int_{\partial\Omega_{\text{in}}} \rho \mathrm{d}\Omega \mathrm{d}t \\ \text{s.t.} \quad & \text{Eq.}\,(1) \\ & G = \int_{\Omega} \alpha \mathrm{d}\Omega - V_{\max}\int_{\Omega} \mathrm{d}\Omega \leq 0 \end{aligned} \tag{10}$$

where $J$ is the objective function. Since the outlet pressure and the mass flow rate are constants in the present case, minimizing the inlet pressure $p_{\text{in}}$ is equivalent to minimizing the pressure drop, as well as equivalent to minimizing the dissipation power [18]. $\partial\Omega_{\text{in}}$ represents the inlet boundary and $W_{\text{in}}$ is its area. $t_{\text{f}}$ is the final moment of the simulation. $G$ is the volume constraint, and $\alpha$ is the design variable, while $V_{\max}$ is the upper limit of the fluid volume. In this work, the design variable is the following binary characteristic function

$$\alpha(\mathbf{x}) = \begin{cases} 0 & \mathbf{x} \text{ is a solid node} \\ 1 & \mathbf{x} \text{ is a fluid node} \end{cases} \tag{11}$$

The dependency of the objective $J$ on the design variable $\alpha$ can be implicitly expressed by expanding the velocity as $\mathbf{u} \leftarrow \alpha\mathbf{u}$ in Eq. (4), leading to the following extended equilibrium distribution function [18]

$$f_i^{\text{eq}} = \omega_i\left\{\rho + \rho_0\left[3(\mathbf{e}_i\cdot\alpha\mathbf{u}) + \frac{9}{2}(\mathbf{e}_i\cdot\alpha\mathbf{u})^2 - \frac{3}{2}|\alpha\mathbf{u}|^2\right]\right\} \tag{12}$$

To make the optimization problem easier to be solved, the above discrete characteristic function is relaxed to a continuous level set function [54]

$$\begin{cases} -1 \leq \phi(\mathbf{x}) < 0 & \mathbf{x} \text{ is a solid node} \\ \phi(\mathbf{x}) = 0 & \mathbf{x} \text{ is a interface node} \\ 0 < \phi(\mathbf{x}) \leq 1 & \mathbf{x} \text{ is a fluid node} \end{cases} \tag{13}$$

During the optimization process, the level set function is updated according to the reaction-diffusion equation [54]

$$\frac{\partial\phi}{\partial\xi} = -K(J' - \sigma\nabla^2\phi) \tag{14}$$

where $\xi$ is the fictitious time, $K$ is the step size and $\sigma$ is the regularization coefficient. $J'$ is the gradient or the sensitivity, defined as the responsible variation of $J$ to an infinite small variation of $\alpha$. In each optimization iteration, the level set function is first updated according to Eq. (14), and then the design variable can be updated according to the following simple mapping

$$\alpha(\mathbf{x}) = \begin{cases} 0 & -1 \le \phi(\mathbf{x}) < 0 \\ 1 & 0 \le \phi(\mathbf{x}) < 1 \end{cases} \tag{15}$$

In the gradient-based TO methods, performing the sensitivity analysis by the adjoint method, or namely the Lagrange multiplier method, is one of the most important tasks. In the adjoint method, the auxiliary adjoint variables are introduced and the Lagrangian function $L$ is constructed by adding the objective function and the inner products of adjoint variables and constraint functions together. By deriving the adjoint equations and solving the adjoint variables, the sensitivity can be conveniently acquired, as to be illustrated later.

*2.3 Continuous adjoint lattice Boltzmann method*

The continuous adjoint method follows the idea of 'first adjoint then discretize', and thus the optimization problem should be formulated in the continuous form, as given in Eq. (10), and the resulting adjoint equation will also be a PDE. Specifically, the Lagrangian function for Eq. (10) will be

$$L = J + \int_t \int_\Omega \sum_{i=0}^{8} f_i^* \left[\frac{\partial f_i}{\partial t} + \mathbf{e}_i \cdot \nabla f_i + \frac{1}{\tau_c}\left(f_i - f_i^{\mathrm{eq}}\right)\right] \mathrm{d}\Omega \mathrm{d}t + \lambda\left(\int_\Omega \alpha \mathrm{d}\Omega - V_{\max} \int_\Omega \mathrm{d}\Omega\right) \tag{16}$$

where $f_i^*$ is the adjoint distribution function for $f_i$, and $\lambda$ is the Lagrange multiplier for volume constraint. To get the adjoint equation, one should take the variational derivative of $L$ over $f_i$, and set this derivative to be zero. In the present case, for arbitrary small variation of $f_i$, i.e., $\delta f_i$, the response variation of $L$ can be expressed by

$$\delta L = \delta J + \int_t \int_\Omega \sum_{i=0}^{8} f_i^* \left[\frac{\partial \delta f_i}{\partial t} + \mathbf{e}_i \cdot \nabla \delta f_i + \frac{1}{\tau_c}\left(\delta f_i - \delta f_i^{\mathrm{eq}}\right)\right] \mathrm{d}\Omega \mathrm{d}t \tag{17}$$

Applying the index perturbation, we have

$$\sum_{i=0}^{8} f_i^* \delta f_i^{\mathrm{eq}} = \sum_{i=0}^{8} \sum_{j=0}^{8} f_i^* \frac{\partial f_i^{\mathrm{eq}}}{\partial f_j} \delta f_j = \sum_{i=0}^{8} \delta f_i \sum_{j=0}^{8} f_j^* \frac{\partial f_j^{\mathrm{eq}}}{\partial f_i} \tag{18}$$

Then using the integral by parts, we have

$$\int_t f_i^* \frac{\partial \delta f_i}{\partial t} \mathrm{d}t = f_i^* \delta f_i |_{(0,t_{\mathrm{f}})} - \int_t \delta f_i \frac{\partial f_i^*}{\partial t} \mathrm{d}t \tag{19}$$

$$\int_\Omega f_i^* \mathbf{e}_i \cdot \nabla \delta f_i \mathrm{d}\Omega = \int_{\partial\Omega} (\mathbf{e}_i \cdot \mathbf{n}) f_i^* \delta f_i \mathrm{d}\Omega - \int_\Omega \delta f_i \mathbf{e}_i \cdot \nabla f_i^* \mathrm{d}\Omega \tag{20}$$

As a result, Eq. (17) can be transformed to

$$\delta L = \delta J + \int_t \int_\Omega \sum_{i=0}^{8} \delta f_i \left[-\frac{\partial f_i^*}{\partial t} - \mathbf{e}_i \cdot \nabla f_i^* + \frac{1}{\tau_c}\left(f_i^* - \sum_{j=0}^{8} f_j^* \frac{\partial f_j^{\mathrm{eq}}}{\partial f_i}\right)\right] \mathrm{d}\Omega \mathrm{d}t + \int_\Omega \sum_{i=0}^{8} f_i^* \delta f_i |_{(0,t_{\mathrm{f}})} \mathrm{d}\Omega + \int_t \int_{\partial\Omega} \sum_{i=0}^{8} (\mathbf{e}_i \cdot \mathbf{n}) f_i^* \delta f_i \mathrm{d}\Omega \mathrm{d}t \tag{21}$$

Considering that

$$\delta J=\delta p_{\mathrm{in}}=\frac{1}{3W_{\mathrm{in}}t_{\mathrm{f}}}\int_t\int_{\partial\Omega_{\mathrm{in}}}\sum_{i=0}^{8}\delta f_i\mathrm{d}\Omega\mathrm{d}t \tag{22}$$

To keep δ$L$ always to be zero, we should have the following adjoint equation satisfied

$$-\frac{\partial f_i^*}{\partial t}-\mathbf{e}_i\cdot\nabla f_i^*+\frac{1}{\tau_{\mathrm{c}}}\left(f_i^*-\sum_{j=0}^{8}f_j^*\frac{\partial f_j^{\mathrm{eq}}}{\partial f_i}\right)=0 \tag{23}$$

with initial condition (at final moment $t$=$t_{\mathrm{f}}$)

$$f_i^*=0 \tag{24}$$

The adjoint boundary conditions can be deduced from the boundary integral terms in Eq. (21). For walls, the adjoint of bounce back is still bounce back, and it is not written here for brevity [18]. For outlet boundary, we should keep

$$\sum_{i=0}^{8}(\mathbf{e}_i\cdot\mathbf{n})f_i^*\delta f_i=0 \tag{25}$$

Substituting the known boundary scheme in Eq. (8) into Eq. (25), we can obtain the adjoint outlet boundary condition as follows

$$\begin{cases}\left(\frac{2}{3}f_2^*+\frac{1}{6}f_5^*+\frac{1}{6}f_6^*\right)\delta f_0=0\\ \left(\frac{2}{3}f_2^*+\frac{2}{3}f_5^*-\frac{1}{3}f_6^*\right)\delta f_1=0\\ \left(\frac{2}{3}f_2^*-\frac{1}{3}f_5^*+\frac{2}{3}f_6^*\right)\delta f_3=0\\ \left(f_4^*-f_2^*+\frac{4}{3}f_2^*+\frac{1}{3}f_5^*+\frac{1}{3}f_6^*\right)\delta f_4=0\\ \left(f_7^*-f_5^*+\frac{4}{3}f_2^*+\frac{1}{3}f_5^*+\frac{1}{3}f_6^*\right)\delta f_7=0\\ \left(f_8^*-f_6^*+\frac{4}{3}f_2^*+\frac{1}{3}f_5^*+\frac{1}{3}f_6^*\right)\delta f_8=0\end{cases} \tag{26}$$

Therefore, the unknown adjoint distribution functions at outlet ($i$=4, 7 and 8) can be calculated by

$$\begin{cases}f_4^*=f_2^*-\left(\frac{4}{3}f_2^*+\frac{1}{3}f_5^*+\frac{1}{3}f_6^*\right)\\ f_7^*=f_5^*-\left(\frac{4}{3}f_2^*+\frac{1}{3}f_5^*+\frac{1}{3}f_6^*\right)\\ f_8^*=f_6^*-\left(\frac{4}{3}f_2^*+\frac{1}{3}f_5^*+\frac{1}{3}f_6^*\right)\end{cases} \tag{27}$$

Similar procedures can be applied to obtain the adjoint inlet boundary condition

$$\begin{cases}\frac{1}{3W_{\mathrm{in}}t_{\mathrm{f}}}\delta f_0=0\\ \left(\frac{1}{3W_{\mathrm{in}}t_{\mathrm{f}}}+\frac{1}{2}f_5^*-\frac{1}{2}f_8^*\right)\delta f_2=0\\ \left(\frac{1}{3W_{\mathrm{in}}t_{\mathrm{f}}}-\frac{1}{2}f_5^*+\frac{1}{2}f_8^*\right)\delta f_4=0\\ \left(\frac{2}{3W_{\mathrm{in}}t_{\mathrm{f}}}-f_1^*+f_3^*\right)\delta f_3=0\\ \left(\frac{2}{3W_{\mathrm{in}}t_{\mathrm{f}}}-f_8^*+f_6^*\right)\delta f_6=0\\ \left(\frac{2}{3W_{\mathrm{in}}t_{\mathrm{f}}}-f_5^*+f_7^*\right)\delta f_7=0\end{cases} \tag{28}$$

Therefore, the unknown adjoint distribution functions at inlet ($i$=3, 6 and 7) can be given by

$$\begin{cases} f_3^* = f_1^* - \frac{2}{3W_{\text{in}}t_{\text{f}}} \\ f_6^* = f_8^* - \frac{2}{3W_{\text{in}}t_{\text{f}}} \\ f_7^* = f_5^* - \frac{2}{3W_{\text{in}}t_{\text{f}}} \end{cases} \tag{29}$$

The above adjoint equation, initial and boundary condition, constitute the complete continuous adjoint lattice Boltzmann model. To solve it, the adjoint equation should be further discretized to

$$f_i^*(\mathbf{x} - \mathbf{e}_i\Delta t, t - \Delta t) - f_i^*(\mathbf{x}, t) = -\frac{1}{\tau}\left(f_i^* - \sum_{j=0}^{8} f_j^* \frac{\partial f_j^{\text{eq}}}{\partial f_i}\right) \tag{30}$$

Then by sequentially performing adjoint collision, adjoint streaming and adjoint boundary condition steps, the adjoint variable $f_i^*$ can be obtained. With the adjoint equation satisfied, namely $\partial L/\partial f_i$=0, the sensitivity can be simplified as

$$J' = \frac{\text{d}L}{\text{d}\alpha} = \frac{\partial L}{\partial \alpha} + \frac{\partial L}{\partial f_i}\frac{\text{d}f_i}{\text{d}\alpha} = \frac{\partial L}{\partial \alpha} \tag{31}$$

Substituting Eq. (12) into Eq. (31), the explicit expression of sensitivity can be given as

$$J' = \frac{\partial L}{\partial \alpha} = \lambda - \frac{1}{\tau}\sum_{i=0}^{8} f_i^* \frac{\partial f_i^{\text{eq}}}{\partial \alpha} = \lambda \underbrace{- \frac{\rho_0}{\tau}\sum_{i=0}^{8} \omega_i f_i^* \left[3(\mathbf{e}_i \cdot \mathbf{u}) + \frac{9}{2}(\mathbf{e}_i \cdot \mathbf{u})^2 - \frac{3}{2}\mathbf{u}^2\right]}_{J_W'} \tag{32}$$

where the Lagrange multiplier $\lambda$ takes the following value in order to accurately enforce the volume constraint [38]

$$\lambda = \frac{\exp(G)}{\int_\Omega \text{d}\Omega}\int_\Omega |J_W' - \sigma\nabla^2\phi| \text{d}\Omega \tag{33}$$

Note that the above continuous adjoint LB model is not rigorously self-consistent if one carefully check the adjoint boundary conditions listed in Eqs. (26) and (28). In Eq. (26), the former three expressions cannot be exactly satisfied unless one assumes that $\delta f_0 = \delta f_1 = \delta f_3 = 0$. Similarly in Eq. (28), the former three equations cannot always be satisfied unless one assumes that $\delta f_0 = \delta f_2 = \delta f_4 = 0$. These assumptions have been utilized in previous studies about the continuous adjoint LBM for open flow systems [39], while it cannot always hold in practice. Therefore, the above adjoint boundary conditions for open systems are intrinsically inconsistent, or singular, or ill-posed. Such singularity has been reported in continuous adjoint models of Navier-Stokes equations [13], and those for DVBEs [45]. The effects of such singularity on the adjoint solution and sensitivity analysis remain unknown, and despite that some assumptions or simplifications are made to the adjoint boundary conditions in the literature to proceed the sensitivity analysis [45, 48], right now it still lacks a rigorous approach to give self-consistent and well-posed adjoint boundary conditions. This work will develop the self-consistent adjoint boundary conditions from the discrete adjoint method, as inspired by the mixed adjoint approach in Ref. [55]. To achieve that, the complete discrete adjoint LB model will be derived and the difference between the adjoint models will be revealed.

*2.4 Discrete adjoint lattice Boltzmann method*

In the discrete adjoint method, the domain and time should be discretized prior to taking the adjoint, thus the governing equation becomes the LBE shown in Eq. (2), and the optimized problem shall be formulated in the following discrete form

$$\begin{aligned} \min_{\boldsymbol{\alpha}} \quad & J = p_{\text{in}} = \frac{1}{3N_{\text{in}}n\Delta t}\sum_{t=0}^{n\Delta t}\sum_{j=1}^{N_{\text{in}}}\rho_j \\ \text{s.t.} \quad & \text{Eq. (2)} \\ & G = \sum_{j=1}^{N}\alpha_j - V_{\max}N \le 0 \end{aligned} \tag{34}$$

where $N_{\text{in}}$ and $N$ are the node numbers of the inlet and the whole domain. Before deriving the adjoint equations, the forward LB model is further illustrated. Eq. (2) is the local evolutional rule for distribution functions. From the global perspective, it can be described by the following individual operations

$$\begin{cases} \text{Collision:} & \mathbf{f}_{\text{C}}(t) = \boldsymbol{\mathcal{C}}[\mathbf{f}(t)] \\ \text{Streaming:} & \mathbf{f}_{\text{S}}(t+\Delta t) = \boldsymbol{\mathcal{S}}[\mathbf{f}_{\text{C}}(t)] \\ \text{Boundary condition:} & \mathbf{f}(t+\Delta t) = \boldsymbol{\mathcal{B}}[\mathbf{f}_{\text{S}}(t+\Delta t)] \end{cases} \tag{35}$$

where $\mathbf{f}$ is the global vector of $f_i$, and $\mathbf{f}_{\text{C}}$ and $\mathbf{f}_{\text{S}}$ are the vector after collision and streaming. Considering that the D2Q9 model is used for 2D pipe bend case, then $\mathbf{f}$, $\mathbf{f}_{\text{C}}$ and $\mathbf{f}_{\text{S}}$ are $9N$ dimensional vectors, in which we put the 9 distribution functions of one point at the neighboring indices. $\boldsymbol{\mathcal{B}}$, $\boldsymbol{\mathcal{S}}$ and $\boldsymbol{\mathcal{C}}$ are the boundary condition, streaming and collision operators for the LBM simulation, and the details are as follows.

(1) Collision operation is purely local but non-linear, and it can be expressed in local form

$$\mathcal{C}_i[\mathbf{f}(\mathbf{x}_j,t)] = f_i(\mathbf{x}_j,t) - \frac{1}{\tau}\left[f_i(\mathbf{x}_j,t) - f_i^{\text{eq}}(\mathbf{x}_j,t)\right] \tag{36}$$

(2) The streaming operation is definitely non-local but strictly linear. Therefore, $\boldsymbol{\mathcal{S}}$ can be expressed by a matrix $\mathbf{S}$, through which we can know that $\mathbf{f}_{\text{S}}(t+\Delta t) = \mathbf{S}\mathbf{f}_{\text{C}}(t)$. To make $\mathbf{S}$ become well-defined, we should make some additional operations. Following Łaniewski-Wołłk and Rokicki [31], we apply periodic boundary treatment when executing streaming. By doing so, the unknown $f_i$ at the boundaries are updated in streaming step. This additional operation will not affect the simulation results at all, since those unknown $f_i$ shall be re-computed by the boundary condition step $\boldsymbol{\mathcal{B}}$. However, this temporary operation make $\mathbf{S}$ become a unitary matrix. One can determine $\mathbf{S}$ according to the explicit understanding to the streaming operation.

In the streaming step, every component in $\mathbf{f}_\mathrm{S}$ comes from and only comes from one certain component in $\mathbf{f}_\mathrm{C}$, and every component in $\mathbf{f}_\mathrm{C}$ goes to and only goes to one certain component in $\mathbf{f}_\mathrm{S}$, therefore, $\mathbf{S}$ is in fact an index perturbation operation, which does not change the elements in $\mathbf{f}_\mathrm{C}$, while only change the element indices of $\mathbf{f}_\mathrm{C}$. Suppose the $i^{\mathrm{th}}$ component in $\mathbf{f}_\mathrm{S}$ comes from the $j^{\mathrm{th}}$ component in $\mathbf{f}_\mathrm{C}$, then $\mathbf{S}_{ij}$=1, and also $\mathbf{S}_{ik}$=0 for $k\neq j$ and $\mathbf{S}_{kj}$=0 for $k\neq i$. Therefore, there will be only 1 element equals to 1 in each row and each column of $\mathbf{S}$. Let us consider

$$(\mathbf{S}\mathbf{S}^{\mathrm{T}})_{ij} = \sum_k \mathbf{S}_{ik}(\mathbf{S}^{\mathrm{T}})_{kj} = \sum_k \mathbf{S}_{ik}\mathbf{S}_{jk} = \begin{cases} 1 & i = j \\ 0 & i \neq j \end{cases} \tag{37}$$

Therefore, we have $\mathbf{S}\mathbf{S}^{\mathrm{T}}$=$\mathbf{I}$, and same for $\mathbf{S}^{\mathrm{T}}\mathbf{S}$=$\mathbf{I}$, $\mathbf{S}$ is unitary.

(3) The Zou-He boundary scheme with purely local boundary treatment is still used, then the boundary condition operation is also purely local but non-linear (affine function), and it can be expressed in local form as shown in Eqs. (7)−(8). For inner nodes, the boundary condition operation is simply identity operation.

The Lagrangian function for Eq. (34) will be

$$L = J + \sum_{t=0}^{n\Delta t} \mathbf{f}_\mathrm{S}^{*\mathrm{T}}(t+\Delta t)\{\mathbf{f}_\mathrm{C}(t) - \boldsymbol{\mathcal{C}}[\mathbf{f}(t)]\} + \sum_{t=0}^{n\Delta t} \mathbf{f}^{*\mathrm{T}}(t+\Delta t)\{\mathbf{f}_\mathrm{S}(t+\Delta t) - \boldsymbol{\mathcal{S}}[\mathbf{f}_\mathrm{C}(t)]\} +$$

$$\sum_{t=0}^{n\Delta t} \mathbf{f}_\mathrm{C}^{*\mathrm{T}}(t+\Delta t)\{\mathbf{f}(t+\Delta t) - \boldsymbol{\mathcal{B}}[\mathbf{f}_\mathrm{S}(t+\Delta t)]\} + \lambda\left(\sum_{j=1}^{N} \alpha_j - V_{\max} N\right) \tag{38}$$

where $n$ is the total number of time steps. $\mathbf{f}^*$ is the vector of the respective adjoint distribution function. To get the adjoint equation, one should calculate the derivatives of $L$ over $\mathbf{f}(t)$, $\mathbf{f}_\mathrm{C}(t)$ and $\mathbf{f}_\mathrm{S}(t)$ individually. For any small variation $\delta\mathbf{f}(t)$, the resulting variation of $L$ will be

$$\delta L = \delta J - \sum_{t=0}^{n\Delta t} \mathbf{f}_\mathrm{S}^{*\mathrm{T}}(t+\Delta t)\delta\boldsymbol{\mathcal{C}}[\mathbf{f}(t)] + \sum_{t=0}^{n\Delta t} \mathbf{f}_\mathrm{C}^{*\mathrm{T}}(t+\Delta t)\delta\mathbf{f}(t+\Delta t) \tag{39}$$

where

$$\delta J = \sum_{t=0}^{n\Delta t} \delta\mathbf{f}^{\mathrm{T}}(t)\left[\frac{\partial J_t}{\partial \mathbf{f}(t)}\right]^{\mathrm{T}} \text{ with } J_t = \frac{1}{3N_{\mathrm{in}} n\Delta t}\sum_{j=1}^{N_{\mathrm{in}}} \rho_j \tag{40}$$

Considering that

$$\sum_{t=0}^{n\Delta t} \mathbf{f}_\mathrm{S}^{*\mathrm{T}}(t+\Delta t)\delta\boldsymbol{\mathcal{C}}[\mathbf{f}(t)] = \sum_{t=0}^{n\Delta t} \mathbf{f}_\mathrm{S}^{*\mathrm{T}}(t+\Delta t)\frac{\partial\boldsymbol{\mathcal{C}}[\mathbf{f}(t)]}{\partial\mathbf{f}(t)}\delta\mathbf{f}(t) =$$

$$\sum_{t=0}^{n\Delta t} \delta\mathbf{f}^{\mathrm{T}}(t)\left\{\frac{\partial\boldsymbol{\mathcal{C}}[\mathbf{f}(t)]}{\partial\mathbf{f}(t)}\right\}^{\mathrm{T}} \mathbf{f}_\mathrm{S}^{*}(t+\Delta t) \tag{41}$$

and

$$\sum_{t=0}^{n\Delta t} \mathbf{f}_\mathrm{C}^{*\mathrm{T}}(t+\Delta t)\delta\mathbf{f}(t+\Delta t) = \sum_{t=1}^{(n+1)\Delta t} \mathbf{f}_\mathrm{C}^{*\mathrm{T}}(t)\delta\mathbf{f}(t) = \sum_{t=1}^{(n+1)\Delta t} \delta\mathbf{f}^{\mathrm{T}}(t)\mathbf{f}_\mathrm{C}^{*}(t) \tag{42}$$

Substituting Eqs. (41)−(42) into Eq. (39) and letting $\delta L$ to be zero, we can get the adjoint collision

$$\mathbf{f}_\mathrm{C}^{*}(t) = \left\{\frac{\partial\boldsymbol{\mathcal{C}}[\mathbf{f}(t)]}{\partial\mathbf{f}(t)}\right\}^{\mathrm{T}} \mathbf{f}_\mathrm{S}^{*}(t+\Delta t) - \left[\frac{\partial J_t}{\partial\mathbf{f}(t)}\right]^{\mathrm{T}} \tag{43}$$

with the adjoint initial condition

$$\mathbf{f}_{\mathrm{C}}^{*}(t)|_{t=(n+1)\Delta t}=0 \tag{44}$$

Then, for any small variation $\delta\mathbf{f}_{\mathrm{C}}(t)$, the resulting variation of $L$ will be

$$\delta L=\sum_{t=0}^{n\Delta t}\mathbf{f}_{\mathrm{S}}^{*\mathrm{T}}(t+\Delta t)\delta\mathbf{f}_{\mathrm{C}}(t)-\sum_{t=0}^{n\Delta t}\mathbf{f}^{*\mathrm{T}}(t+\Delta t)\delta\boldsymbol{\mathcal{S}}[\mathbf{f}_{\mathrm{C}}(t)] \tag{45}$$

Note that

$$\sum_{t=0}^{n\Delta t}\mathbf{f}_{\mathrm{S}}^{*\mathrm{T}}(t+\Delta t)\delta\mathbf{f}_{\mathrm{C}}(t)=\sum_{t=0}^{n\Delta t}\delta\mathbf{f}_{\mathrm{C}}^{\mathrm{T}}(t)\mathbf{f}_{\mathrm{S}}^{*}(t+\Delta t) \tag{46}$$

and

$$\sum_{t=0}^{n\Delta t}\mathbf{f}^{*\mathrm{T}}(t+\Delta t)\delta\boldsymbol{\mathcal{S}}[\mathbf{f}_{\mathrm{C}}(t)]=\sum_{t=0}^{n\Delta t}\delta\mathbf{f}_{\mathrm{C}}^{\mathrm{T}}(t)\mathbf{S}^{\mathrm{T}}\mathbf{f}^{*}(t+\Delta t) \tag{47}$$

Substituting Eqs. (46)–(47) into Eq. (45) and letting $\delta L$ to be zero, we can get the adjoint streaming

$$\mathbf{f}_{\mathrm{S}}^{*}(t+\Delta t)=\mathbf{S}^{\mathrm{T}}\mathbf{f}^{*}(t+\Delta t) \tag{48}$$

Finally, for any small variation $\delta\mathbf{f}_{\mathrm{S}}(t+\Delta t)$, the resulting variation of $L$ will be

$$\delta L=\sum_{t=0}^{n\Delta t}\mathbf{f}^{*\mathrm{T}}(t+\Delta t)\delta\mathbf{f}_{\mathrm{S}}(t+\Delta t)-\sum_{t=0}^{n\Delta t}\mathbf{f}_{\mathrm{C}}^{*\mathrm{T}}(t+\Delta t)\delta\boldsymbol{\mathcal{B}}[\mathbf{f}_{\mathrm{S}}(t+\Delta t)] \tag{49}$$

Considering that

$$\sum_{t=0}^{n\Delta t}\mathbf{f}^{*\mathrm{T}}(t+\Delta t)\delta\mathbf{f}_{\mathrm{S}}(t+\Delta t)=\sum_{t=0}^{n\Delta t}\delta\mathbf{f}_{\mathrm{S}}^{\mathrm{T}}(t+\Delta t)\mathbf{f}^{*}(t+\Delta t) \tag{50}$$

and

$$\sum_{t=0}^{n\Delta t}\mathbf{f}_{\mathrm{C}}^{*\mathrm{T}}(t+\Delta t)\delta\boldsymbol{\mathcal{B}}[\mathbf{f}_{\mathrm{S}}(t+\Delta t)]=\sum_{t=0}^{n\Delta t}\mathbf{f}_{\mathrm{C}}^{*\mathrm{T}}(t+\Delta t)\frac{\partial\boldsymbol{\mathcal{B}}[\mathbf{f}_{\mathrm{S}}(t+\Delta t)]}{\partial\mathbf{f}_{\mathrm{S}}(t+\Delta t)}\delta\mathbf{f}_{\mathrm{S}}(t+\Delta t)=\sum_{t=0}^{n\Delta t}\delta\mathbf{f}_{\mathrm{S}}^{\mathrm{T}}(t+\Delta t)\left\{\frac{\partial\boldsymbol{\mathcal{B}}[\mathbf{f}_{\mathrm{S}}(t+\Delta t)]}{\partial\mathbf{f}_{\mathrm{S}}(t+\Delta t)}\right\}^{\mathrm{T}}\mathbf{f}_{\mathrm{C}}^{*}(t+\Delta t) \tag{51}$$

Substituting Eqs. (50)–(51) into Eq. (49) and letting $\delta L$ to be zero, we can get the adjoint boundary condition

$$\mathbf{f}^{*}(t+\Delta t)=\left\{\frac{\partial\boldsymbol{\mathcal{B}}[\mathbf{f}_{\mathrm{S}}(t+\Delta t)]}{\partial\mathbf{f}_{\mathrm{S}}(t+\Delta t)}\right\}^{\mathrm{T}}\mathbf{f}_{\mathrm{C}}^{*}(t+\Delta t) \tag{52}$$

Summarizing the above formulations, we can determine the solution strategy for the discrete adjoint LBM as follows.

(1) Initially $t=(n+1)\Delta t$, at that time we have

$$\mathbf{f}_{\mathrm{C}}^{*}(t)|_{t=(n+1)\Delta t}=0 \tag{53}$$

(2) Then we execute adjoint boundary condition

$$\mathbf{f}^{*}(t+\Delta t)=\left\{\frac{\partial\boldsymbol{\mathcal{B}}[\mathbf{f}_{\mathrm{S}}(t+\Delta t)]}{\partial\mathbf{f}_{\mathrm{S}}(t+\Delta t)}\right\}^{\mathrm{T}}\mathbf{f}_{\mathrm{C}}^{*}(t+\Delta t) \tag{54}$$

(3) Then the adjoint streaming

$$\mathbf{f}_{\mathrm{S}}^{*}(t+\Delta t)=\mathbf{S}^{\mathrm{T}}\mathbf{f}^{*}(t+\Delta t) \tag{55}$$

(4) Finally the adjoint collision

$$\mathbf{f}_{\mathrm{C}}^{*}(t)=\left\{\frac{\partial\boldsymbol{\mathcal{C}}[\mathbf{f}(t)]}{\partial\mathbf{f}(t)}\right\}^{\mathrm{T}}\mathbf{f}_{\mathrm{S}}^{*}(t+\Delta t)-\left[\frac{\partial J_t}{\partial\mathbf{f}(t)}\right]^{\mathrm{T}} \tag{56}$$

Note that the above formulations are in the global perspective, to render the adjoint solution as efficient as the forward solution, now we should recover the adjoint solution steps to their local forms.

For adjoint boundary condition, let us denote

$$\mathbf{B} = \frac{\partial \boldsymbol{\mathcal{B}}[\mathbf{f}_{\mathrm{S}}(t+\Delta t)]}{\partial \mathbf{f}_{\mathrm{S}}(t+\Delta t)} \tag{57}$$

Since $\boldsymbol{\mathcal{B}}$ only involves local operation, **B** will be block diagonal matrix, and its transpose, $\mathbf{B}^{\mathrm{T}}$ will also be a block diagonal matrix. As a result, the adjoint boundary treatment is also purely local operation. In an inner point (of course, most of the points will be inner points), the local blocks will be simply identity matrix, and the adjoint is also identity, thus $f_i^* = f_{\mathrm{C},i}^*$. In a local boundary point $j$=1, 2, …, $N$, e.g., inlet, outlet or wall, we can derive the expression of **B** by hand. For example, at inlet we have

$$\mathbf{B}_j = \begin{bmatrix} 1 & 0 & 0 & 0 & 0 & 0 & 0 & 0 & 0 \\ 0 & 0 & 0 & 1 & 0 & 0 & 0 & 0 & 0 \\ 0 & 0 & 1 & 0 & 0 & 0 & 0 & 0 & 0 \\ 0 & 0 & 0 & 1 & 0 & 0 & 0 & 0 & 0 \\ 0 & 0 & 0 & 0 & 1 & 0 & 0 & 0 & 0 \\ 0 & 0 & -0.5 & 0 & 0.5 & 0 & 0 & 1 & 0 \\ 0 & 0 & 0 & 0 & 0 & 0 & 1 & 0 & 0 \\ 0 & 0 & 0 & 0 & 0 & 0 & 0 & 1 & 0 \\ 0 & 0 & 0.5 & 0 & -0.5 & 0 & 1 & 0 & 0 \end{bmatrix} \tag{58}$$

Substituting $\mathbf{B}_j$ into Eq. (54), we can get the adjoint inlet boundary condition

$$\begin{aligned} & f_1^* = f_5^* = f_8^* = 0 \\ & f_0^* = f_{\mathrm{C},0}^* \\ & f_2^* = f_{\mathrm{C},2}^* - \tfrac{1}{2} f_{\mathrm{C},5}^* + \tfrac{1}{2} f_{\mathrm{C},8}^* \\ & f_4^* = f_{\mathrm{C},4}^* + \tfrac{1}{2} f_{\mathrm{C},5}^* - \tfrac{1}{2} f_{\mathrm{C},8}^* \\ & f_3^* = f_{\mathrm{C},3}^* + f_{\mathrm{C},1}^* \\ & f_6^* = f_{\mathrm{C},6}^* + f_{\mathrm{C},8}^* \\ & f_7^* = f_{\mathrm{C},7}^* + f_{\mathrm{C},5}^* \end{aligned} \tag{59}$$

Similarly we can reach the adjoint outlet boundary condition

$$\begin{aligned} & f_2^* = f_5^* = f_6^* = 0 \\ & f_0^* = f_{\mathrm{C},0}^* - \tfrac{2}{3} f_{\mathrm{C},2}^* - \tfrac{1}{6} f_{\mathrm{C},5}^* - \tfrac{1}{6} f_{\mathrm{C},6}^* \\ & f_1^* = f_{\mathrm{C},1}^* - \tfrac{2}{3} f_{\mathrm{C},2}^* - \tfrac{2}{3} f_{\mathrm{C},5}^* + \tfrac{1}{3} f_{\mathrm{C},6}^* \\ & f_3^* = f_{\mathrm{C},3}^* - \tfrac{2}{3} f_{\mathrm{C},2}^* + \tfrac{1}{3} f_{\mathrm{C},5}^* - \tfrac{2}{3} f_{\mathrm{C},6}^* \\ & f_4^* = f_{\mathrm{C},4}^* - \tfrac{1}{3} f_{\mathrm{C},2}^* - \tfrac{1}{3} f_{\mathrm{C},5}^* - \tfrac{1}{3} f_{\mathrm{C},6}^* \\ & f_7^* = f_{\mathrm{C},7}^* - \tfrac{4}{3} f_{\mathrm{C},2}^* + \tfrac{2}{3} f_{\mathrm{C},5}^* - \tfrac{1}{3} f_{\mathrm{C},6}^* \\ & f_8^* = f_{\mathrm{C},8}^* - \tfrac{4}{3} f_{\mathrm{C},2}^* - \tfrac{1}{3} f_{\mathrm{C},5}^* + \tfrac{2}{3} f_{\mathrm{C},6}^* \end{aligned} \tag{60}$$

For walls, $\mathbf{B}_j$ is also an index perturbation operation, and it is a symmetric matrix. Thus the adjoint for bounce back is still bounce back, and the expression will be omitted here.

For adjoint streaming, note that $\mathbf{S}$ is unitary, and $\mathbf{S}^{\mathrm{T}}=\mathbf{S}^{-1}$, so the adjoint streaming is the reversed streaming, which can be expressed in the local form as follows

$$f_{\mathrm{S},i}^{*}\left(\mathbf{x}_j,t+\Delta t\right)=f_i^{*}\left(\mathbf{x}_j+\mathbf{e}_i\Delta t,t+\Delta t\right) \tag{61}$$

For adjoint collision, let us denote

$$\mathbf{C}=\frac{\partial \boldsymbol{\mathcal{C}}[\mathbf{f}(t)]}{\partial \mathbf{f}(t)} \tag{62}$$

Since $\boldsymbol{\mathcal{C}}$ only involves local operation, $\mathbf{C}$ will be block diagonal matrix, and its transpose, $\mathbf{C}^{\mathrm{T}}$ will also be a block diagonal matrix. As a result, the adjoint collision is also purely local operation. For any inner point $j$=1, 2, …, $N$, the local matrix operator is

$$\mathbf{C}_j=\frac{\partial \boldsymbol{\mathcal{C}}[\mathbf{f}(\mathbf{x}_j,t)]}{\partial \mathbf{f}(\mathbf{x}_j,t)}=\frac{\partial\left\{\mathbf{f}(\mathbf{x}_j,t)-\frac{1}{\tau}[\mathbf{f}(\mathbf{x}_j,t)-\mathbf{f}^{\mathrm{eq}}(\mathbf{x}_j,t)]\right\}}{\partial \mathbf{f}(\mathbf{x}_j,t)}=\mathbf{I}-\frac{1}{\tau}\left[\mathbf{I}-\frac{\partial \mathbf{f}^{\mathrm{eq}}(\mathbf{x}_j,t)}{\partial \mathbf{f}(\mathbf{x}_j,t)}\right] \tag{63}$$

$\mathbf{C}_j$ is a block diagonal matrix (each block in 9×9 dimension). Hence the adjoint collision at point $j$ will be

$$\mathbf{f}_{\mathrm{C}}^{*}\left(\mathbf{x}_j,t\right)=\mathbf{C}_j{}^{\mathrm{T}}\mathbf{f}_{\mathrm{S}}^{*}\left(\mathbf{x}_j,t+\Delta t\right)-\left[\frac{\partial J_t}{\partial \mathbf{f}(\mathbf{x}_j,t)}\right]^{\mathrm{T}}=\mathbf{f}_{\mathrm{S}}^{*}\left(\mathbf{x}_j,t+\Delta t\right)-\frac{1}{\tau}\left\{\mathbf{f}_{\mathrm{S}}^{*}\left(\mathbf{x}_j,t+\Delta t\right)-\left[\frac{\partial \mathbf{f}^{\mathrm{eq}}(\mathbf{x}_j,t)}{\partial \mathbf{f}(\mathbf{x}_j,t)}\right]^{\mathrm{T}}\mathbf{f}_{\mathrm{S}}^{*}\left(\mathbf{x}_j,t+\Delta t\right)\right\}-\left[\frac{\partial J_t}{\partial \mathbf{f}(\mathbf{x}_j,t)}\right]^{\mathrm{T}} \tag{64}$$

The adjoint collision in any discrete direction $i$ will be

$$f_{\mathrm{C},i}^{*}\left(\mathbf{x}_j,t\right)=f_{\mathrm{S},i}^{*}\left(\mathbf{x}_j,t+\Delta t\right)-\frac{1}{\tau}\left\{f_{\mathrm{S},i}^{*}\left(\mathbf{x}_j,t+\Delta t\right)-\sum_{m=0}^{8}\frac{\partial f_m^{\mathrm{eq}}(\mathbf{x}_j,t)}{\partial f_i(\mathbf{x}_j,t)}f_{\mathrm{S},m}^{*}\left(\mathbf{x}_j,t+\Delta t\right)\right\}-\frac{\partial J_t}{\partial f_i(\mathbf{x}_j,t)} \tag{65}$$

It can be found that the solution strategy of the discrete adjoint LBM is also very similar to the primal LBM, that is, by repeatedly performing the local adjoint boundary condition, the non-local adjoint streaming, and the local adjoint collision, the adjoint variable $\mathbf{f}^{*}$ can be solved. The sensitivity in the present method can be given as

$$J'\left(\mathbf{x}_j\right)=\frac{\partial L}{\partial \alpha_j}=-\mathbf{f}_{\mathrm{S}}^{*\mathrm{T}}\left(\mathbf{x}_j\right)\frac{\partial\{\boldsymbol{\mathcal{C}}[\mathbf{f}(\mathbf{x}_j)]\}}{\partial \alpha}+\lambda=-\frac{1}{\tau}\sum_{i=0}^{8}f_{\mathrm{S},i}^{*}\left(\mathbf{x}_j\right)\frac{\partial f_i^{\mathrm{eq}}(\mathbf{x}_j)}{\partial \alpha}+\lambda=\underbrace{-\frac{\rho_0}{\tau}\sum_{i=0}^{8}\omega_i f_{\mathrm{S},i}^{*}\left(\mathbf{x}_j\right)\left[3\left(\mathbf{e}_i\cdot\mathbf{u}\left(\mathbf{x}_j\right)\right)+\frac{9}{2}\left(\mathbf{e}_i\cdot\mathbf{u}\left(\mathbf{x}_j\right)\right)^2-\frac{3}{2}\mathbf{u}^2\left(\mathbf{x}_j\right)\right]}_{J_W'(\mathbf{x}_j)}+\lambda \tag{66}$$

where the Lagrange multiplier $\lambda$ should be expressed in the discretized manner now

$$\lambda=\frac{\exp(G)}{N}\sum_{j=1}^{N}\left|J_W'(\mathbf{x}_j)-\sigma\nabla^2\phi(\mathbf{x}_j)\right| \tag{67}$$

One can note that the discrete adjoint method leads to a fully-consistent adjoint model, meaning that no self-contradiction appears and no assumption is needed during the derivation.

The benefits brought by such self-consistency will be discussed and demonstrated later.

### *2.5 Comparison between two adjoint models*

Comparing the adjoint models derived from the two adjoint methods, one can note that the most significant difference between them is the self-consistency pertaining to the adjoint boundary conditions, that is, the continuous adjoint method gives inconsistent adjoint boundary conditions while the discrete adjoint method can obtain fully consistent adjoint model. Specifically, the differences in the adjoint boundary conditions are the following points.

**Point 1:** The resulting expressions of the adjoint boundary conditions are different, as can be known by comparing Eqs. (29) and (59), or Eqs. (27) and (60).

**Point 2:** The effects of objective functions are different. In the continuous adjoint LBM, the related derivative term is added into the boundary condition, as shown in Eq. (29); but in the discrete adjoint LBM, the related derivative term becomes a source term added in the boundary node, as shown in Eq. (65).

**Point 3:** The execution orders of adjoint collision, adjoint streaming and adjoint boundary conditions are different. In the continuous adjoint LBM, the order is first collision, then streaming and finally boundary condition, just the same order as the primal LBM. However, in the discrete adjoint method, the execution order becomes first boundary condition, then streaming and finally collision, completely reversed to the primal LBM, as a direct result of the global transpose operation.

Except the above differences, the discrete adjoint model is exactly the same as its continuous counterpart. For example, the adjoint collisions given in Eqs. (30) and (65) are exactly the same if the adjoint source term in Eq. (65) in excluded, and the adjoint streaming shown in Eqs. (30) and (61) is actually the same revered streaming. Besides, the sensitivity expressions presented in Eqs. (32) and (66) are exactly the same. Therefore, it can be asserted that the only difference between the two adjoint methods lies in the different adjoint boundary conditions they generated. The expression, the effect of objective function and the execution order, which are the essential factors related to the adjoint boundary conditions, are totally different, corresponding to totally different consistency performance of the two adjoint methods.

## 3 Results and discussion

Now we have figured out the difference between the two adjoint models, it is natural to further ask that will such difference bring different numerical performance when performing the sensitivity analysis. In fact, the primal motivation of the present research is to improve the numerical stability of the continuous adjoint LBM. It has been reported in the literature [43], as well as been testified in early times by the authors that the above continuous adjoint LBM has a rather poor numerical stability, much inferior than the primal LBM simulation. Specifically, the adjoint solution gets diverged once *Re* approaches nearly 25 [43], while the primal LBM simulation can give converged solution at a much higher *Re*, e.g., 200. However, when the continuous adjoint method is applied to perform sensitivity analysis in closed systems, such as in the natural convection case, the resulting adjoint LBM shows excellent numerical stability at high Grashof number (*Gr*), and delicate optimized heat sinks at high *Gr* can be given, as shown in Fig. 3 [2, 14]. The difference between open and closed systems lies in the boundary conditions. In the open system, the Zou-He boundary scheme leads to inconsistent adjoint boundary conditions at inlet and outlet, while in the closed system, all boundaries are no-slip walls implemented by bounce-back scheme, which can lead to naturally consistent adjoint boundary condition [14, 18]. Therefore, it is reasonable to infer that the poor numerical stability of the continuous adjoint LBM is rooted in the inconsistency of the continuous adjoint boundary conditions at the open boundaries. Furthermore, it can be inferred or expected that the discrete adjoint method with fully consistent adjoint boundary conditions can improve the numerical stability of the adjoint solution, which will be demonstrated by the 2D and 3D pipe bend cases in the below.

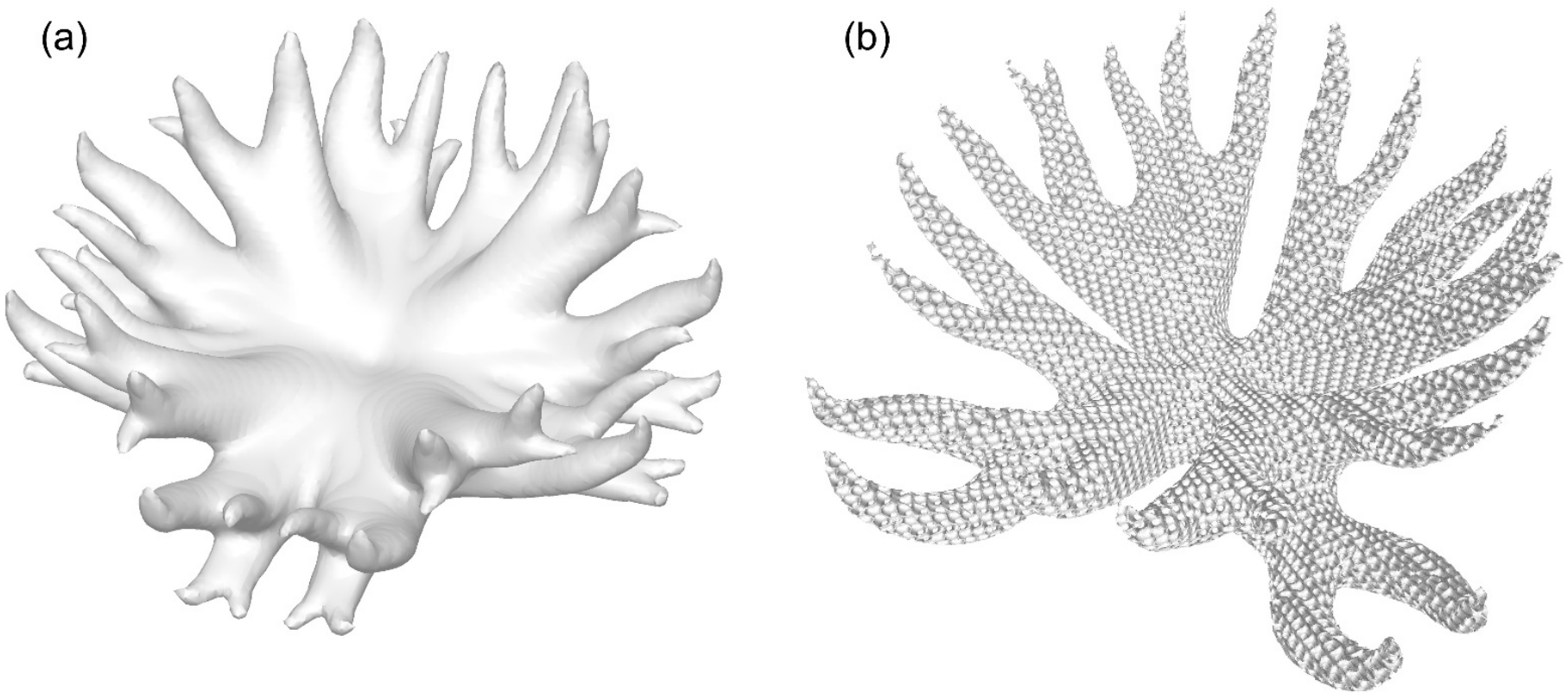

**Figure 3.** Optimized heat sinks at high *Gr*. (a) Solid heat sink at $Gr$=1.6×10$^6$ [2]; (b) porous heat sink at $Gr$=1.2×10$^5$ [14].

*3.1 2D pipe bend*

According to the numerical methods introduced in Section 2, the pipe bend designs under different *Re* can be obtained. Fig. 4 presents the convergence history of the objective function (in lattice unit, lu), while the optimized pipe bend designs is shown in Fig. 5. The initial geometry is the empty chamber, and the maximum fluid volume fraction is set as 0.08$\pi$. When the optimization starts, the objective value first increases since more solid are gradually added into the design domain. After the volume constraint is satisfied, the pipe shape is further subtly modified to lower the objective value. The case with low *Re* as 0.2 is the same benchmark case presented by Borrvall and Petersson [5], and their results are also presented in Fig. 5 to make a direct comparison. It can be found that the current results given by both continuous and discrete adjoint LBM agree well with those of Borrvall and Petersson [5]. When *Re* is increased, the flow channel becomes more curved, and similar trend has been reported in Ref. [38]. Such change of flow channel morphology reflects the change of dominant pressure loss. When *Re* is low, the pressure loss from the viscous friction dominates, the thus the flow channel should be as straight as possible to shorten the flow path and to reduce the viscous loss. Yet when *Re* is higher, the inertia of the influent fluid is stronger and the inertial loss becomes dominant, thus the flow channel should be more curved to gradually deflect the flow and reduce the inertial loss.

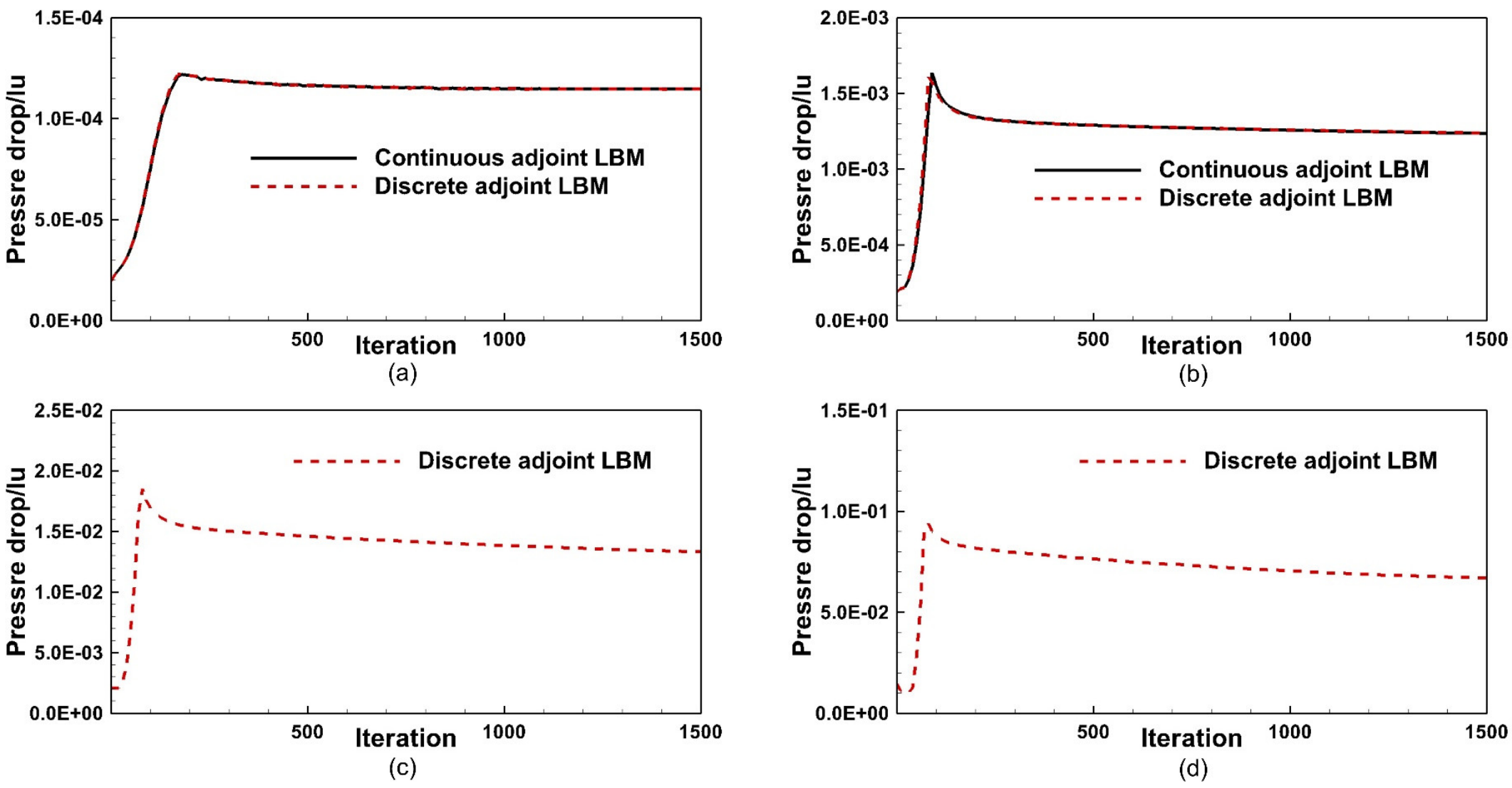

**Figure 4.** Convergence of the objective value in the 2D pipe bend case. (a) *Re*=0.2; (b) *Re*=2; (c) *Re*=20; (d) *Re*=200.

It can be noted that the continuous adjoint LBM cannot give optimized design when *Re* exceeds 20 due to the divergence of the adjoint solution, while the discrete adjoint LBM can get reasonable optimized designs, showing a much superior numerical stability. To show this more clearly, a stability test is performed, where the maximum velocity that can keep the program stable under certain relaxation time is recorded. Here the initial geometry of the pipe bend case is adopted to perform the calculation, and the test is deemed stable if it does not diverge within 1,000,000 iterations [56]. The results are presented in Fig. 6. It can be confirmed that the numerical stability of the discrete adjoint LBM is much superior to that of the continuous adjoint LBM, with a nearly 10 times higher maximum velocity achieved under the same relaxation time. Meanwhile, it can be noted that the discrete adjoint LBM is just as stable as the primal LBM, which is also the theoretically highest stability that can be achieved by the adjoint solution. Note that some high velocities could be unusable for practical simulation due to the compressibility error, and they are used here for stability test only [56]. As mentioned in Section 2, the only difference between the continuous adjoint LBM and the discrete adjoint LBM is the implementation of adjoint boundary conditions, in which the continuous adjoint boundary conditions are intrinsically inconsistent while the discrete ones are completely self-consistent. Therefore, it can be confirmed that the poor numerical stability of the continuous adjoint LBM is rooted in the inconsistency of the continuous adjoint boundary conditions at the open boundaries, and the discrete adjoint method with fully consistent adjoint boundary conditions can greatly improve the numerical stability of the adjoint solution, rendering it as high as the primal LBM.

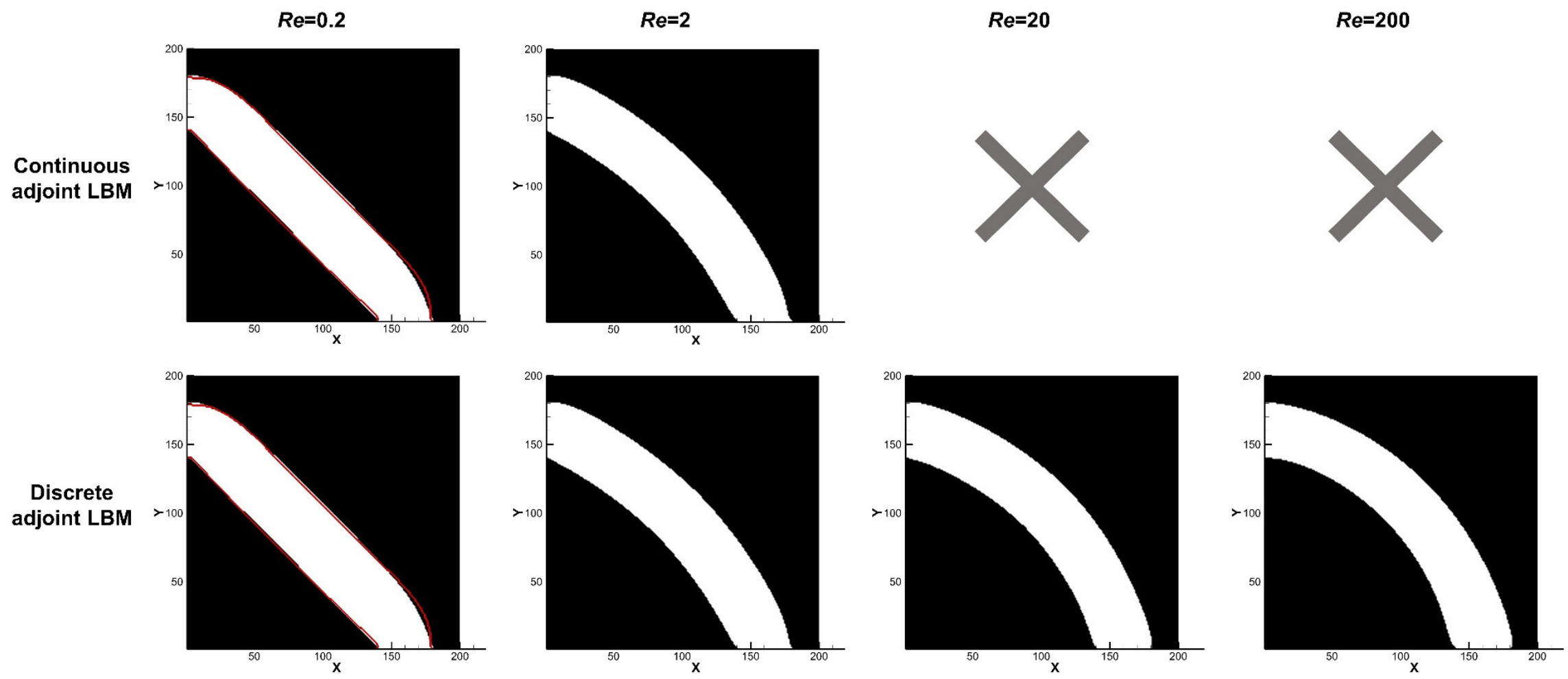

**Figure 5.** Optimized 2D pipe bends at various *Re* using continuous and discrete adjoint LBM. Red lines represent the results from Borrvall and Petersson [5].

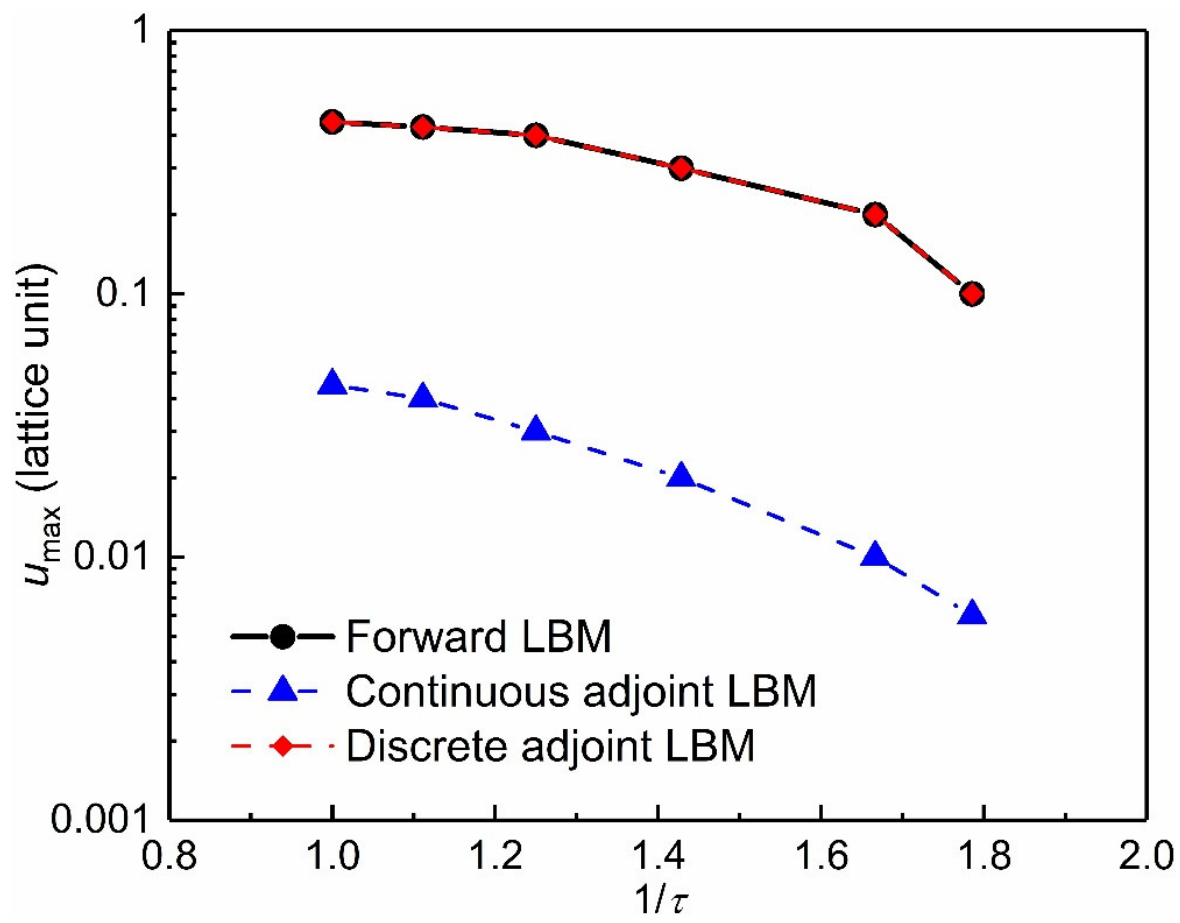

**Figure 6.** Stability test to the 2D continuous and discrete adjoint LBM.

Besides the numerical stability, the accuracy of sensitivity analysis and the computational efficiency of the two adjoint methods are further studied. The sensitivity distributions given by the two adjoint methods at *Re*=0.2 and 2 are displayed in Fig. 7, and the line distributions along the diagonal are extracted to compare with the results given by the FDM. It can be seen that the sensitivity predicted by the continuous adjoint LBM deviates from that of the FDM, and such deviation becomes more significant when *Re* gets higher. In contrast, the discrete adjoint LBM always give accurate sensitivity results. Therefore, it can be concluded that the inconsistency of continuous adjoint LBM causes some errors in the sensitivity analysis, and the fully consistent discrete adjoint LBM is able to eliminate such errors. The computation time of the two adjoint methods at *Re*=0.2 and 2 is also recorded in Table 1, in which all of the optimizations are performed in the same high-performance cluster with 64 AMD EPYC 7452 CPU cores. It can be observed that the discrete adjoint LBM spends similar or little less computation time compared to the continuous adjoint LBM. It is worth pointing out that the two adjoint methods are expected to have the same computational efficiency from the theoretical perspective, since their solution strategies are almost the same and their main computational expenses are similar. For instance, the collision step that contributes most of the local computation burden, and the streaming step that is responsible for most of the communication expenses, are exactly the same for the two adjoint methods, and only the boundary condition step involving little amount of computation and communication tasks is different for the two adjoint methods. Therefore, it can be summarized that the present discrete LBM is as efficient as the continuous adjoint LBM. Fig. 8 further compares the convergence behavior of the forward LBM solution and the two adjoint LBM solutions, where the residual of quantity $I$ (e.g. average velocity magnitude), which is defined as $|I(t)-I(t-1000\Delta t)|/|I(t)|$, is plotted. It can be found that the two adjoint LBM solutions have similar convergence behavior to the forward LBM solution, and the discrete adjoint LBM can give similar convergence speed compared with the continuous adjoint LBM, leading to similar computation time of the two adjoint methods.

**Table 1**. Computation time of the 2D pipe bend case

| *Re* | Computation time | |
|---|---|---|
| | Continuous adjoint LBM | Discrete adjoint LBM |
| 0.2 | 760 s | 648 s |
| 2 | 585 s | 440 s |

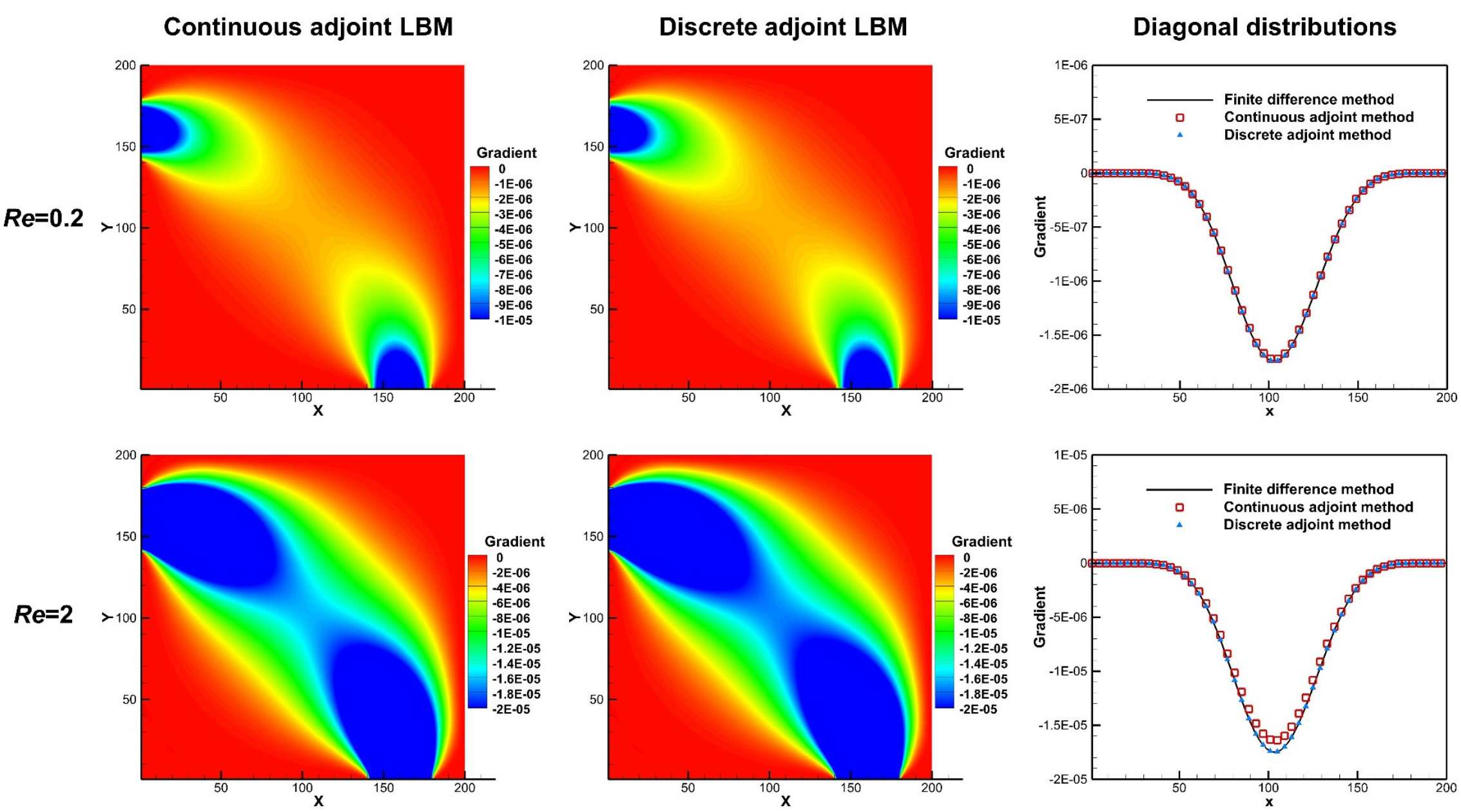

**Figure 7.** Comparison of sensitivities given by the 2D continuous and discrete adjoint LBM.

*3.2 3D pipe bend*

To further confirm the above conclusions in 3D circumstances, the 3D pipe bend case illustrated in Fig. 9 is further investigated. Here the formulation of the optimization problem is still in the form of Eq. (34), with $V_{\max}$ being reset to $0.004\pi^2$. The D3Q19 stencil shown in Fig. 2 is applied here, and the Zou-He boundary scheme is utilized to deal with open boundaries. The 3D adjoint models given by the two adjoint methods are basically the same as their 2D counterpart, while the related adjoint boundary conditions can be derived in a similar manner, and the results are put in Supplementary material A considering their verbosity. The convergence history of the objective value is plotted in Fig. 10, and the optimized designs under different *Re* are plotted in Fig. 11. Similar to the 2D cases, the objective value first increases and then decreases until convergence is reached. The benchmark case with *Re*=0.2 is exactly the same case in Challis and Guest [57], and one can note the similar straight channel is obtained here, proving its validity. Again similar to the 2D cases, the channels get bended when *Re* is increased, and the continuous adjoint LBM gets diverged when *Re* reaches 20. The stability test shown in Fig. 12 clearly proves the much better numerical stability of the discrete adjoint LBM, confirming the conclusions given by the 2D case. Fig. 13 compares the sensitivity distributions obtained from the two adjoint methods, where the discrete adjoint LBM again outperforms the continuous one regarding the accuracy of sensitivity. Table 2 records the computation time of the two methods using 128 AMD EPYC 7452 CPU cores, which shows that the computational efficiencies of the two adjoint methods are basically the same. Fig. 14 again compares the convergence behavior of the forward and two adjoint LBM solutions, and similar convergence speed of the two adjoint methods can be confirmed. Therefore, the above test results from 3D cases well support the conclusions given by the 2D cases, and it can be concluded that for open flow systems, the fully-consistent discrete adjoint LBM is able to improve the numerical stability and accuracy compared to the inconsistent continuous adjoint LBM, while at the same time it keeps the high computational efficiency.

**Table 2**. Computation time of the 3D pipe bend case

| *Re* | Computation time | |
|---|---|---|
| | Continuous adjoint LBM | Discrete adjoint LBM |
| 0.2 | 52 min | 52 min |
| 2 | 35 min | 35 min |

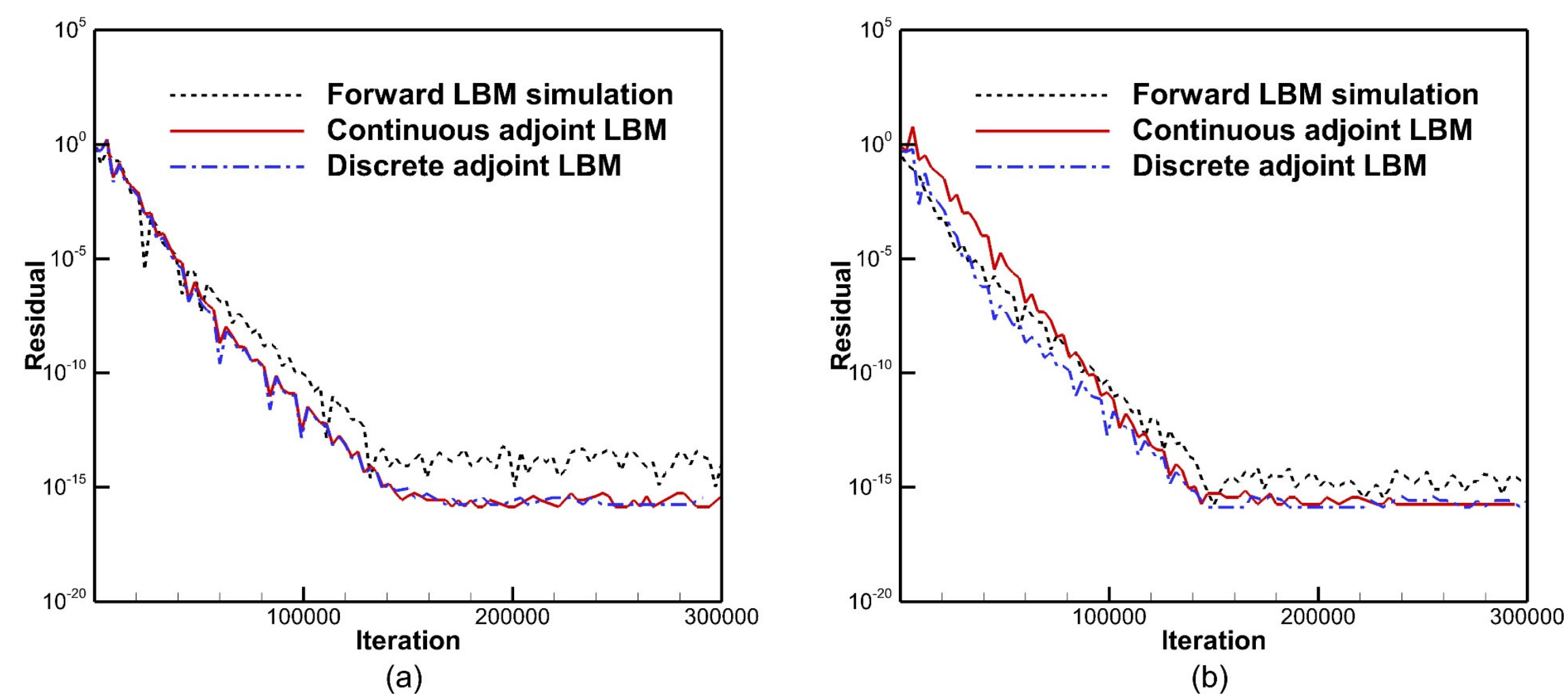

**Figure 8.** Convergence of the forward and adjoint solutions at the first optimization step in the 2D pipe bend case. (a) *Re*=0.2; (b) *Re*=2.

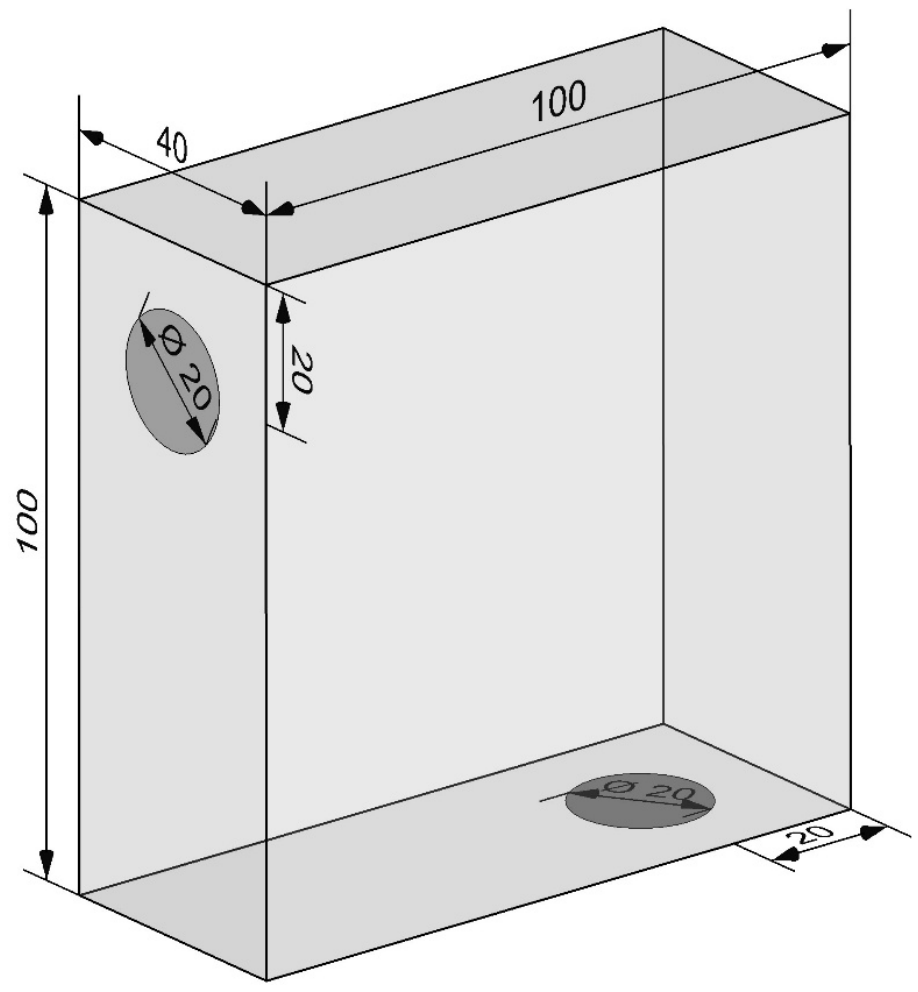

**Figure 9.** Schematic of the 3D pipe bend case.

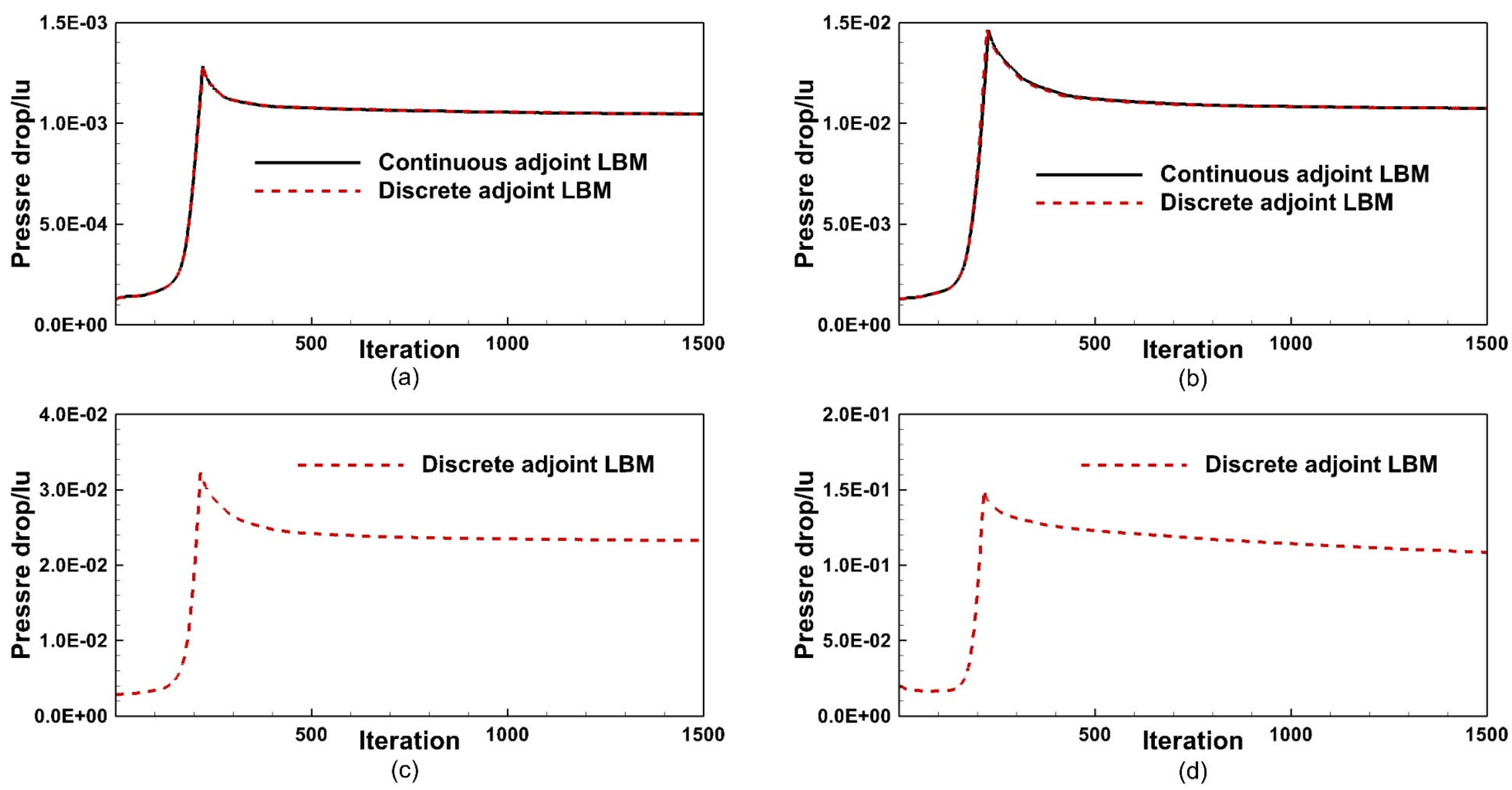

**Figure 10.** Convergence of the objective value in the 3D pipe bend case. (a) *Re*=0.2; (b) *Re*=2; (c) *Re*=20; (d) *Re*=80.

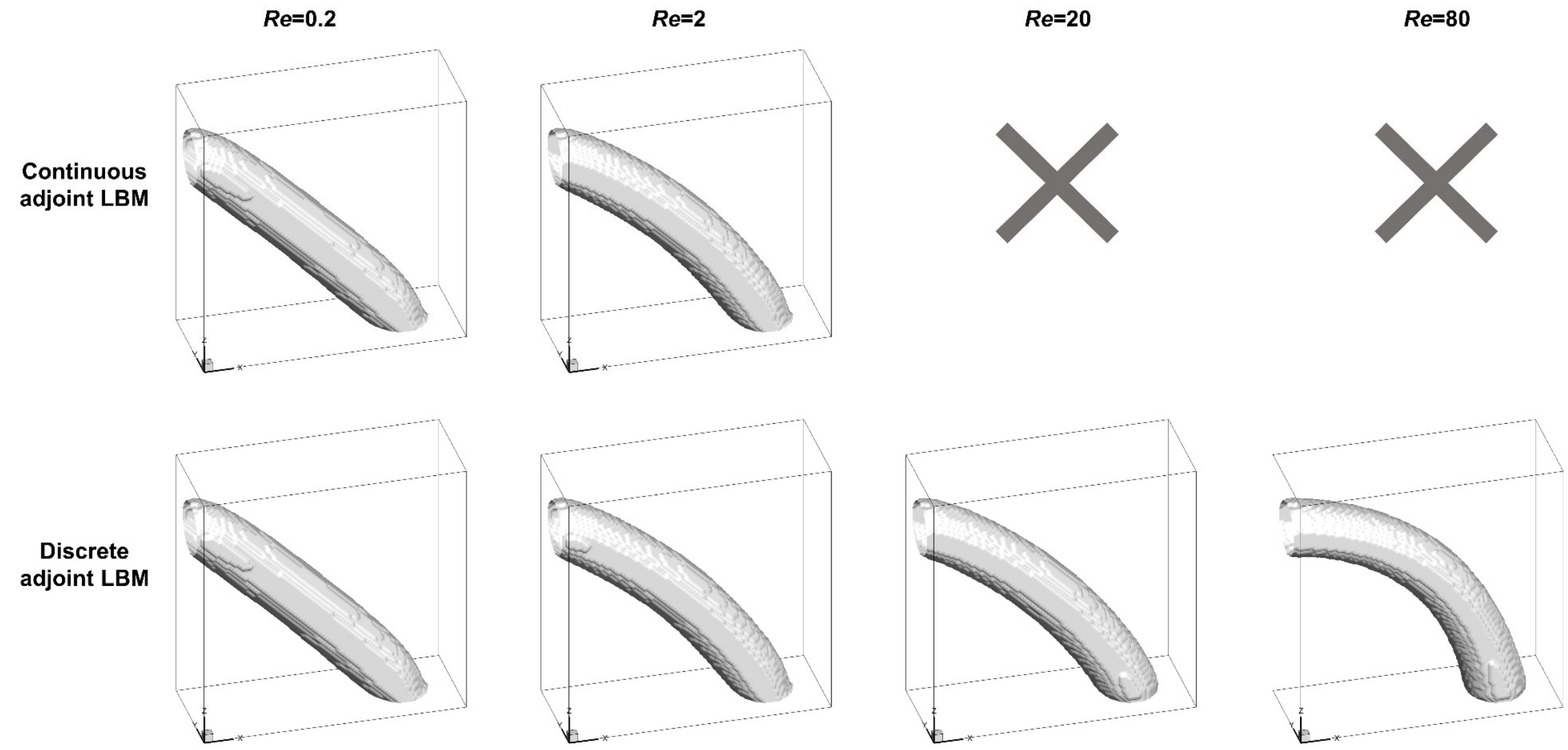

**Figure 11.** Optimized 3D pipe bends at various *Re* using continuous and discrete adjoint LBM.

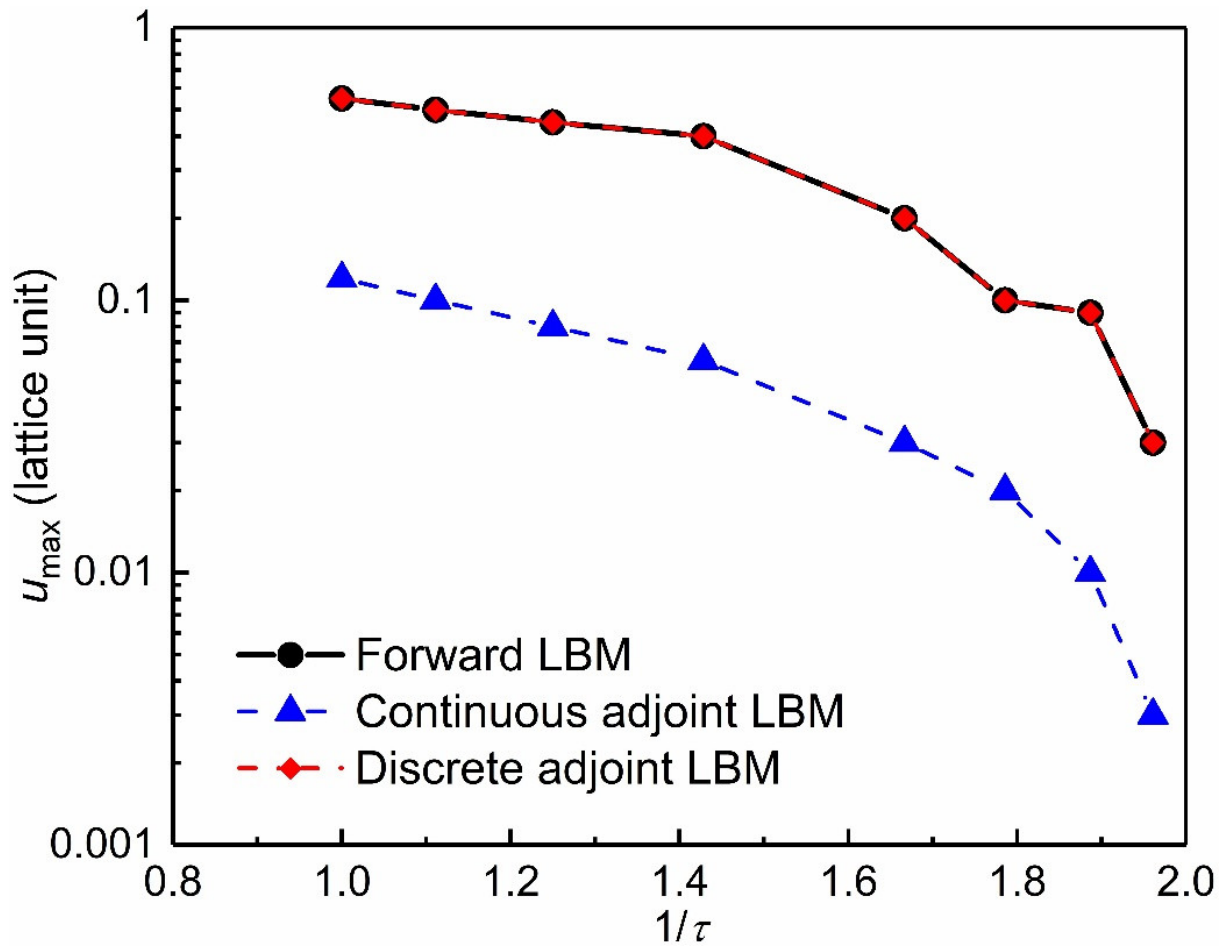

**Figure 12.** Stability test to the 3D continuous and discrete adjoint LBM.

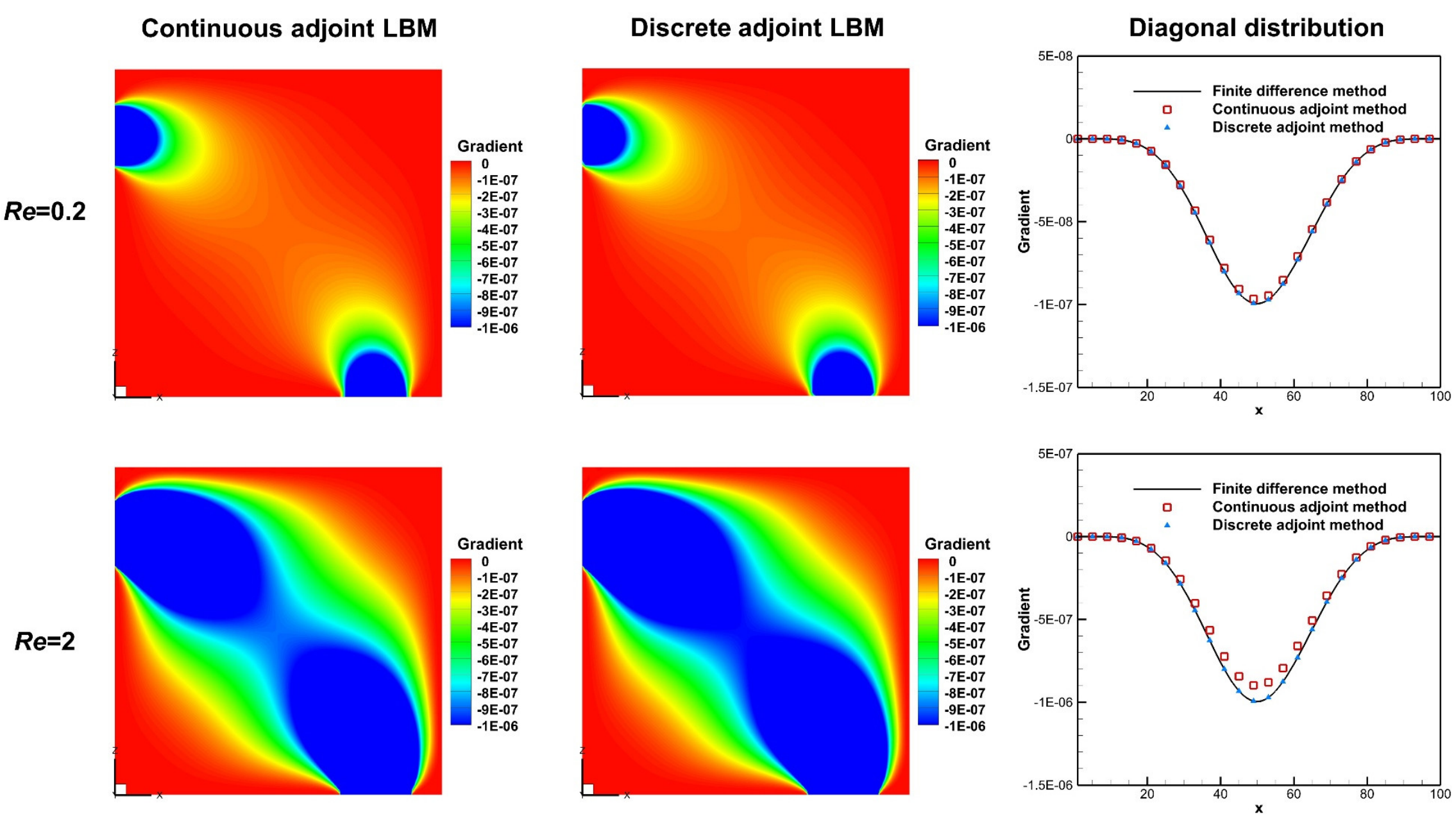

**Figure 13.** Comparison of sensitivities given by the 3D continuous and discrete adjoint LBM.

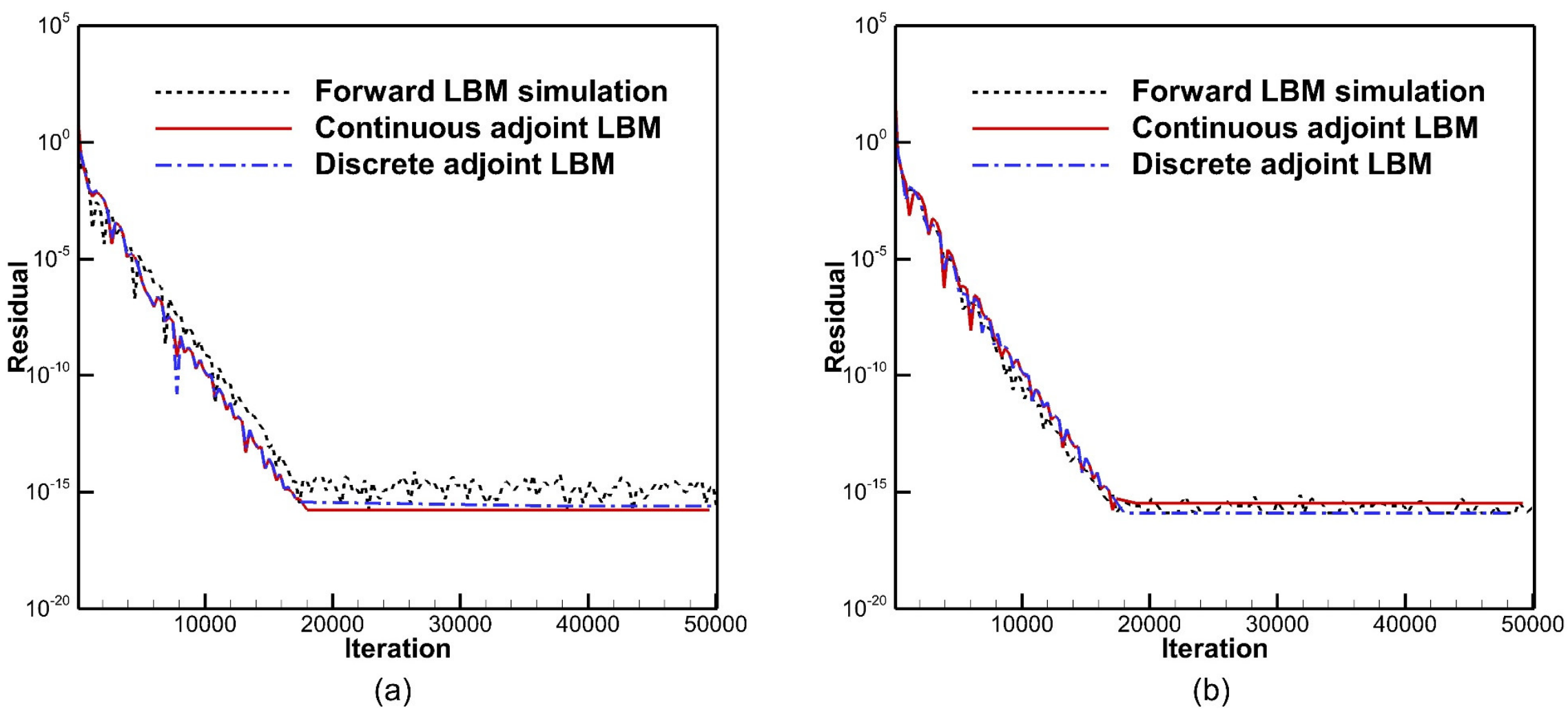

**Figure 14.** Convergence of the forward and adjoint solutions at the first optimization step in the 3D pipe bend case. (a) *Re*=0.2; (b) *Re*=2.

### *3.3 3D microchannel heat sink*

We have demonstrated that the discrete adjoint LBM can solve the numerical instability issues of the continuous adjoint LBM by deriving the fully consistent adjoint boundary conditions at open boundaries. One of the most important applications for deriving such robust approach for sensitivity analysis is to design the active microchannel heat sinks using TO, which has been the hotspot in the research community for a long time [58]. To demonstrate the versatility of the discrete adjoint LBM developed here, this section presents the design of microchannel heat sink under different *Re*.

The schematic of the present case in displayed in Fig. 15, which is similar to the case presented by Yaji et al. [58]. The cold fluid with high pressure $p_1$ and low temperature $T_0$ is injected into the design domain (in green) to evacuate the heat that is generated and conducted inside the solid, and it flows out from the exit with low pressure $p_0$ and thermal fully developed status. Other boundaries are all no-slip and adiabatic walls. The steady-state laminar fluid flow and heat transfer are considered, and the dimensionless governing equations are as follows [58]

$$\begin{gathered} \nabla' \cdot \mathbf{u}' = 0 \\ (\mathbf{u}' \cdot \nabla')\mathbf{u}' = -\nabla' p' + \frac{1}{Re}\nabla'^2 \mathbf{u}' \\ RePr(\mathbf{u}' \cdot \nabla')T' = \nabla'^2 T' + \beta'(\alpha)(1 - T') \end{gathered} \tag{68}$$

where the variables with prime are the dimensionless variables, for example, $T'$ is the dimensionless temperature. $Pr$ is the Prandtl number, and $\beta'$ is the dimensionless heat generation coefficient, which can be expressed as [58]

$$\beta'(\alpha) = \beta'_{\max}(1 - \alpha) \tag{69}$$

where $\beta'_{\max}$ is the maximum heat generation coefficient in the solid. Besides, considering the symmetry, only a quarter of the domain, as colored in Fig. 15, is utilized to perform the simulation, and symmetrical boundary condition is applied at the cutting planes. Since a fixed pressure instead of a given velocity is assigned at the inlet, *Re* is computed by a characteristic velocity estimated from the pressure drop [58], which can be written as

$$Re = \frac{\overline{U}L}{\nu} \quad \text{with} \quad \overline{U} = \sqrt{\frac{\Delta p}{\rho_0}} = \sqrt{\frac{p_1 - p_0}{\rho_0}} \tag{70}$$

To resolve the convective heat transfer problem, another distribution function $g_i$ should be introduced, whose evolution function can be given by

$$g_i(\mathbf{x} + \mathbf{e}_i \Delta t, t + \Delta t) - g_i(\mathbf{x}, t) = -\frac{1}{\tau_g}\left[g_i(\mathbf{x}, t) - g_i^{\mathrm{eq}}(\mathbf{x}, t)\right] + \omega_i Q(\mathbf{x}, t) \tag{71}$$

where $\tau_g$ is the relaxation time and $g_i^{\mathrm{eq}}$ is the equilibrium distribution function. $Q$ is the heat source. Since the heat transfer LBM used here is the same one presented in our previous work [14], it is omitted here for brevity and one can refer to Ref. [14] for more details. The objective of this design problem is to maximize the heat evacuated by the cold fluid, and considering the energy conservation, it is equivalent to maximize the total heat generation. Hence the discrete form of the optimization problem writes [58]

$$\begin{aligned} \min_{\boldsymbol{\alpha}} \quad & J = -\textstyle\sum_{j=1}^{N} Q(\mathbf{x}_j) = -\sum_{j=1}^{N} \beta'(\alpha_j)(1 - T'_j) \\ \text{s.t.} \quad & \text{Eqs.}\,(2) \ \text{and} \ (71) \\ & G = \textstyle\sum_{j=1}^{N} \alpha_j - V_{\max} N \le 0 \end{aligned} \tag{72}$$

Similarly, the Lagrange multiplier method introduced in Section 2 should be applied to conduct the sensitivity analysis. Here the derivation process is omitted for brevity, and the resulting adjoint models and the sensitivity expression are put in Supplementary material B.

In the present case, *Pr* is set to 7.1 since the working fluid is water. $\beta'_{max}$ is set to $1\times10^{-5}$ and the thermal conductivity ratio between solid and fluid is 100. *Re* is varied by changing the inlet pressure, and three values of *Re*, namely *Re*=60, 120 and 240, are considered. Note that the continuous adjoint LBM will diverge when *Re* is larger than ~20, thus the discrete adjoint LBM is absolutely necessary here.

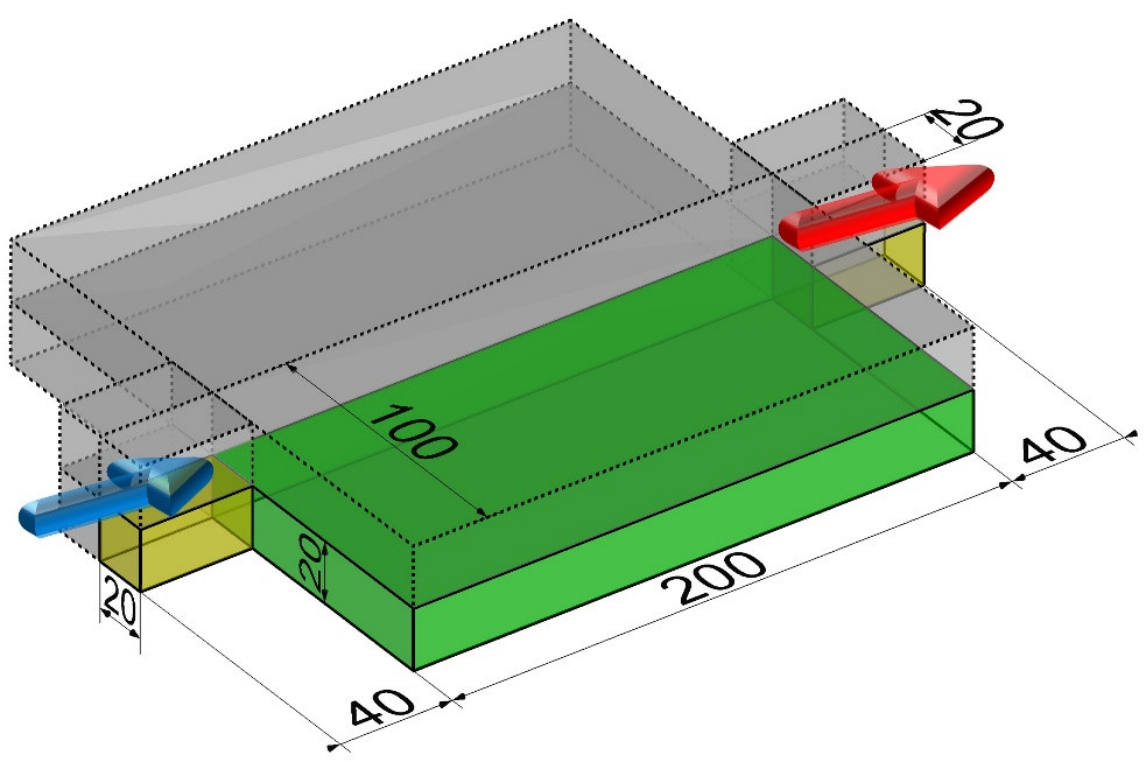

**Figure 15.** Schematic of the 3D microchannel heat sink case.

3.3.1 Results of the base case

The convergence history of the objective values is plotted in Fig. 16, and the optimized designs and the corresponding velocity and temperature distributions are displayed in Fig. 17. The topology optimization starts from the initial structure with full fluid in the design domain, and the solid material is gradually introduced as the optimization proceeds, during which the objective value may either increase or decrease. After the volume constraint gets satisfied, the design is further optimized and the objective value (i.e., amount of heat transfer) increases until convergence is reached. It can be found that the flow channels become much more complicated when increasing *Re*, which is qualitatively similar to the conclusion in Ref. [58]. Increasing *Re* also leads to larger flow velocity and enhanced heat transfer or equivalently lower temperature in the solid, as can be confirmed from Fig. 17. When *Re*=60, the flow channel is simple and only a large solid part exists in the central domain. Such simple and wide flow channels help to reduce the flow resistance and allow larger flow flux under a fixed pressure drop, which can facilitate the heat exchange. When *Re* gets larger, the flow channel becomes narrower and more winding, and meanwhile the solid parts become smaller and more dispersed. Such variation of optimized design can be explained by the following fact. When *Re* increases, the flow rate increases and more fluid is able to flow through the narrow and winding flow channels, where the narrow channels lead to increased heat exchange surface area and the winding channels correspond to a longer flow path. Meanwhile, the more dispersed solid parts in the domain can distribute the heat source more uniformly in the domain, leading to a greatly enhanced heat transfer.

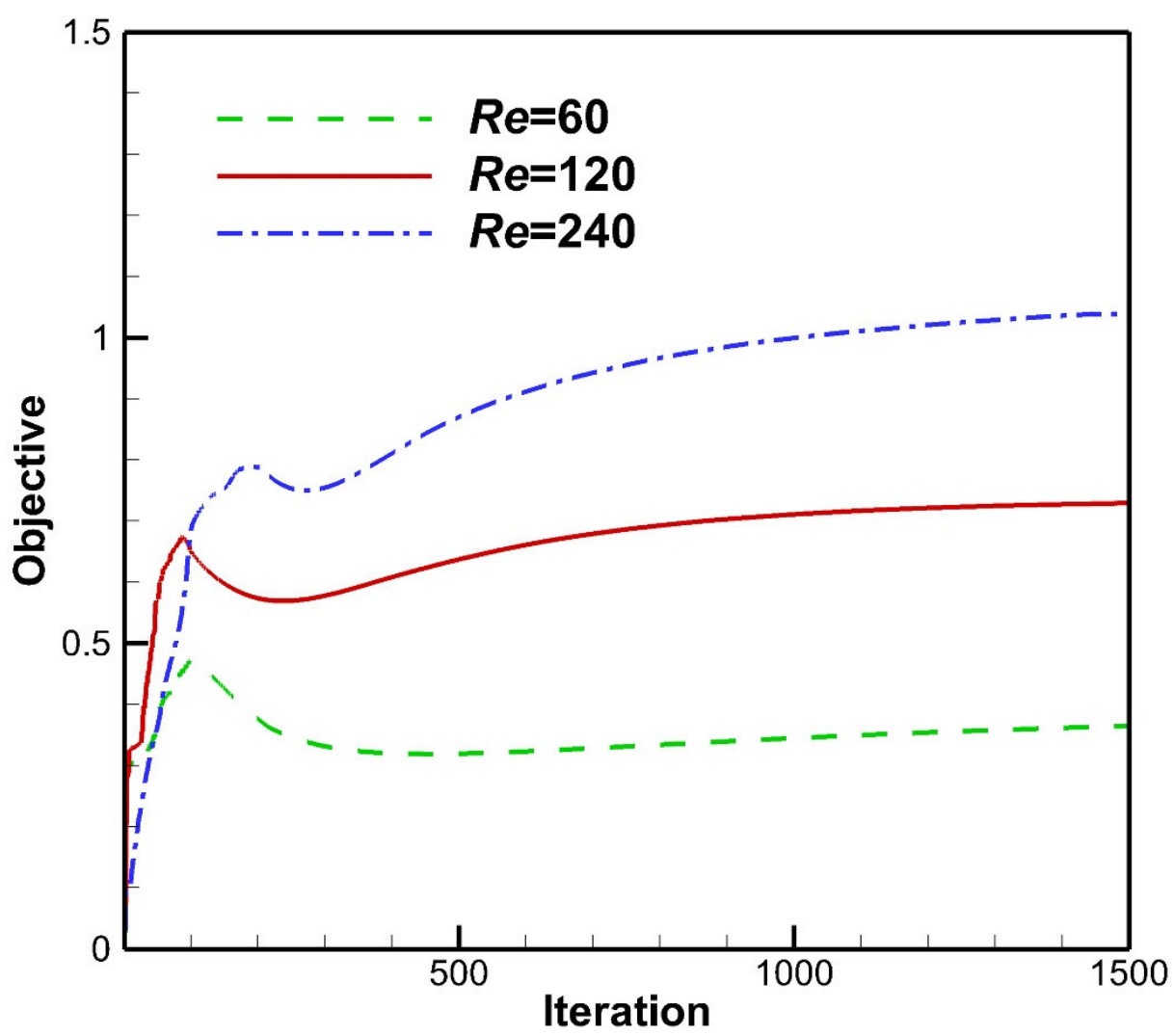

**Figure 16.** Convergence of the objective values in the 3D heat sink cases.

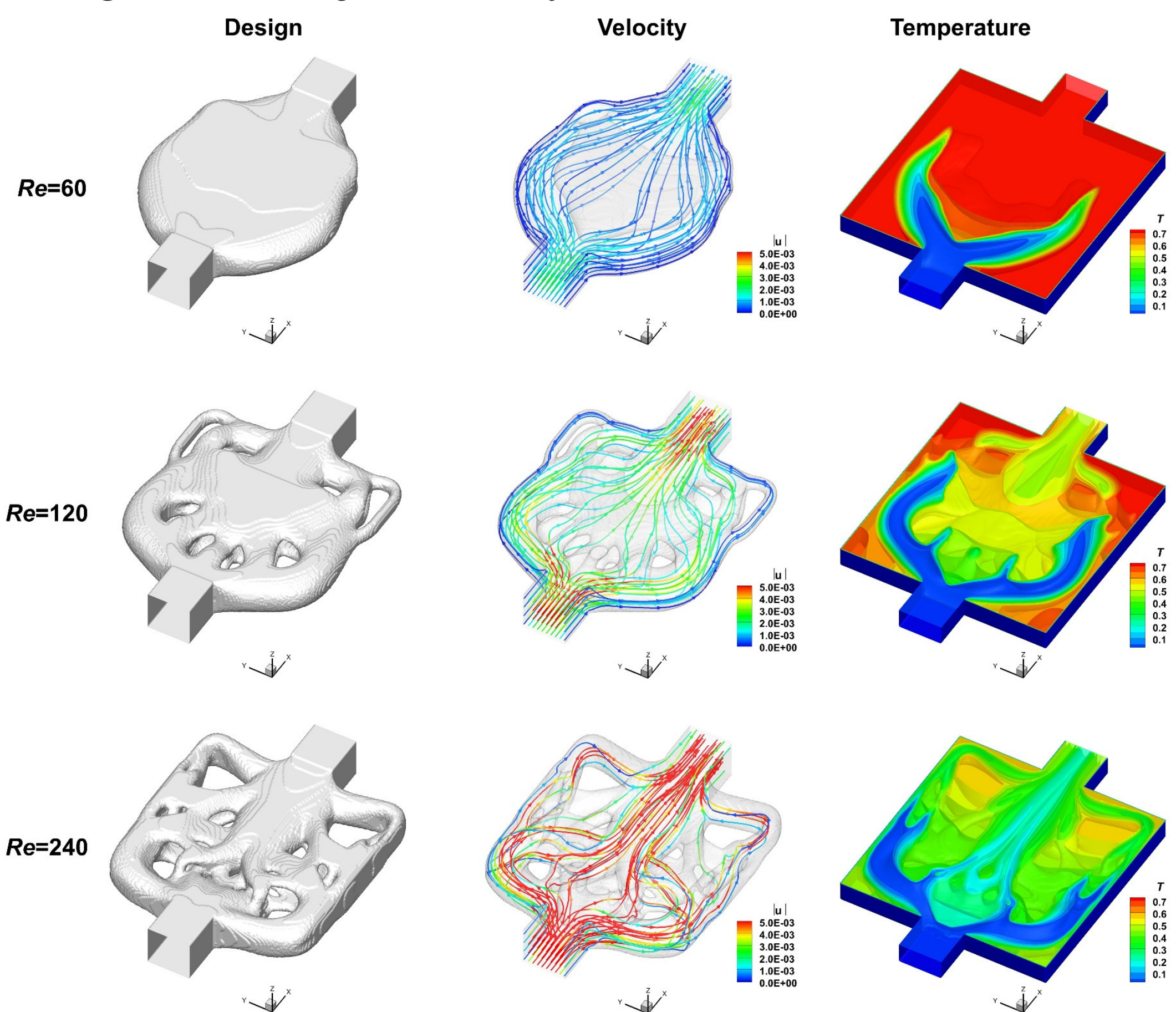

**Figure 17.** Optimized 3D microchannel heat sinks of the base case and the respective velocity and temperature distributions.

To prove the effectiveness of the optimized designs and validity of the above analysis, a cross check to the optimized designs is performed and results are listed in Table 3. It can be observed that the design optimized under certain *Re* has a better performance than other designs for its specific *Re*, well supporting the analysis above. In summary, the above results demonstrate the necessity and versatility of the discrete adjoint LBM developed in the present work.

It is also worth mentioning that the present study focuses on filling the gap of numerical stability between the adjoint solution and forward solution. Because the simplest single-relaxation-time collision model with quite limited stability is applied as the forward solver here, the values of *Re* involved in the present work cannot be too high. However, since the present work renders the adjoint solution as stable as the forward solution, the numerical stability can be readily improved and much higher *Re* can be achieved when some more stable collision models (e.g. the multiple-relaxation-time collision model) is incorporated [2], while it is out of the scope of this work and can be left to future research.

**Table 3**. Cross check of the optimized designs under different *Re*

| *Re* | Objective value | | |
|---|---|---|---|
| | Optimized design at *Re*=60 | Optimized design at *Re*=120 | Optimized design at *Re*=240 |
| 60 | **0.372** | 0.329 | 0.178 |
| 120 | 0.693 | **0.733** | 0.528 |
| 240 | 0.888 | 1.015 | **1.042** |

### 3.3.2 Effects of the thermal conductivity ratio

Effects of the thermal conductivity ratio are studied by changing its value from 100 to 1. The heat sink structures optimized under the three characteristic *Re* values are plotted in Fig. 18. It can be noted that the optimized design becomes more complex when *Re* increases, similar to the findings in the base case. Since the solid thermal conductivity is much lower than that of the base case, one can find that the outlet temperature of the present case is obviously lower, mainly due to the weakened heat transfer brought by poorer heat conduction in the solid. For the same reason, the temperature variation or the temperature gradient in the solid parts becomes much more obvious. The dimensionless objective values of the present case are 0.306, 0.598 and 0.852 respectively, much lower than those of the base case, indicating a deteriorated heat transfer. It is also interesting to note that when *Re*=60, the optimized design of the present case is quite similar to that in the base case, because the low *Re* prevents the further complication of the design. However, when *Re*=240, the design under low thermal conductivity ratio is more complex than that of the base case, because dispersing solid parts in the design domain is more preferable when the solid thermal conductivity is low, and meanwhile the high *Re* allows more complex flow channels.

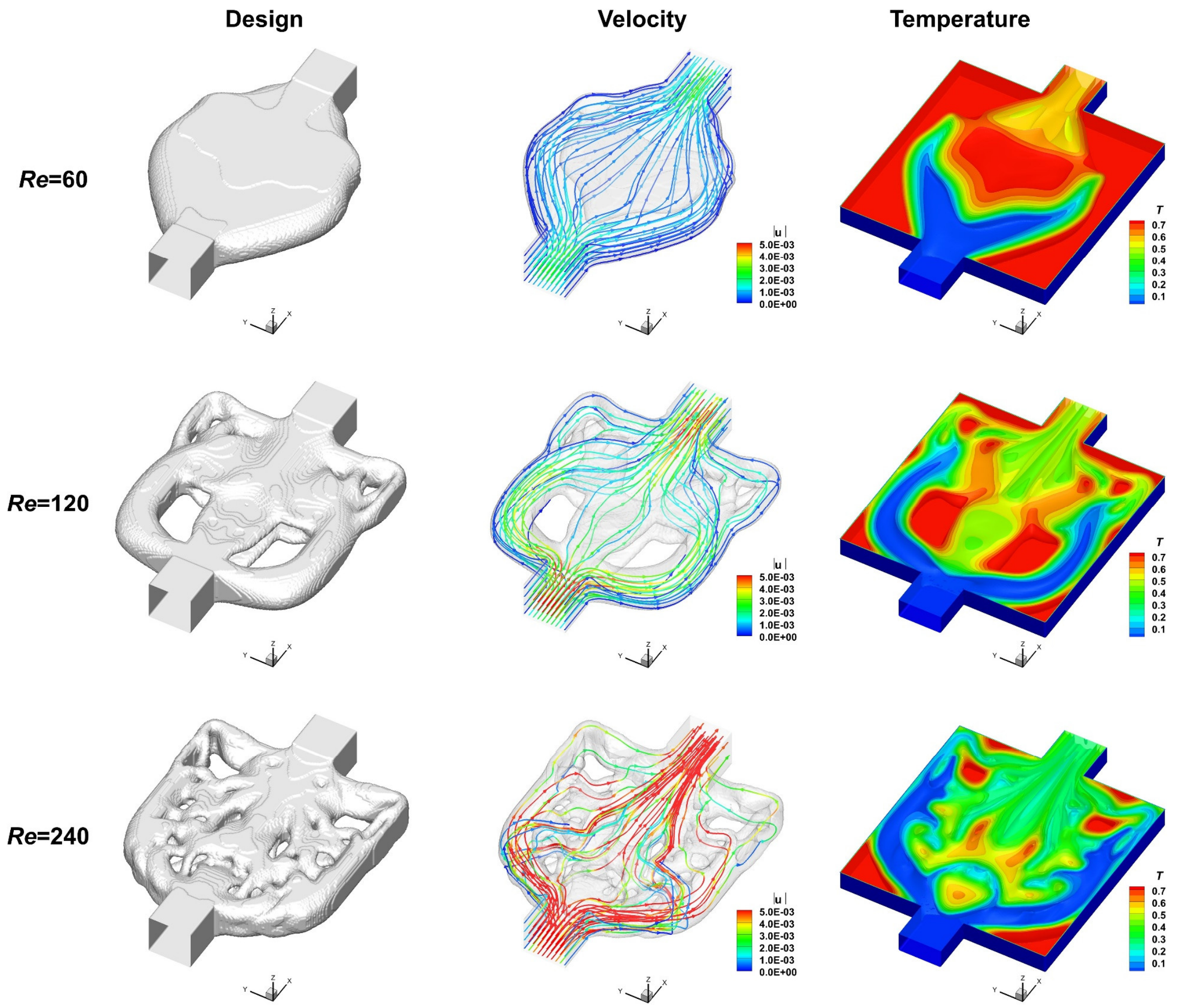

**Figure 18.** Optimized 3D microchannel heat sinks under low thermal conductivity ratio and the respective velocity and temperature distributions.

### 3.3.3 Effects of the fluid volume fraction

In the present case, the fluid volume fraction is set to 0.3 while it is 0.5 in the base case. Fig. 19 displays the optimized designs under three characteristic *Re* values. It can be observed that the flow channels become simpler due to the decrease of fluid volume. Specifically, the main flow channels get narrower and some complex secondary flow channels disappear. Still the flow channels are more winding and dispersed in the domain along with the increase of *Re*, which can lengthen the flow path for heat exchange, facilitating the overall heat transfer. In this case, the dimensionless objective values of the present case are 0.176, 0.540 and 0.796 respectively, much lower than those of the base case, indicating a deteriorated heat transfer. Although more solid that could work as heat source is available in the present case, the total heat transfer amount decreases since less fluid can be injected into the domain to evacuate heat when the pressure drop keeps unchanged while the flow channels get narrower. In other words, the complex coupling between the fluid flow and heat transfer processes makes the topology optimization of forced convection heat transfer more intricate to analyze.

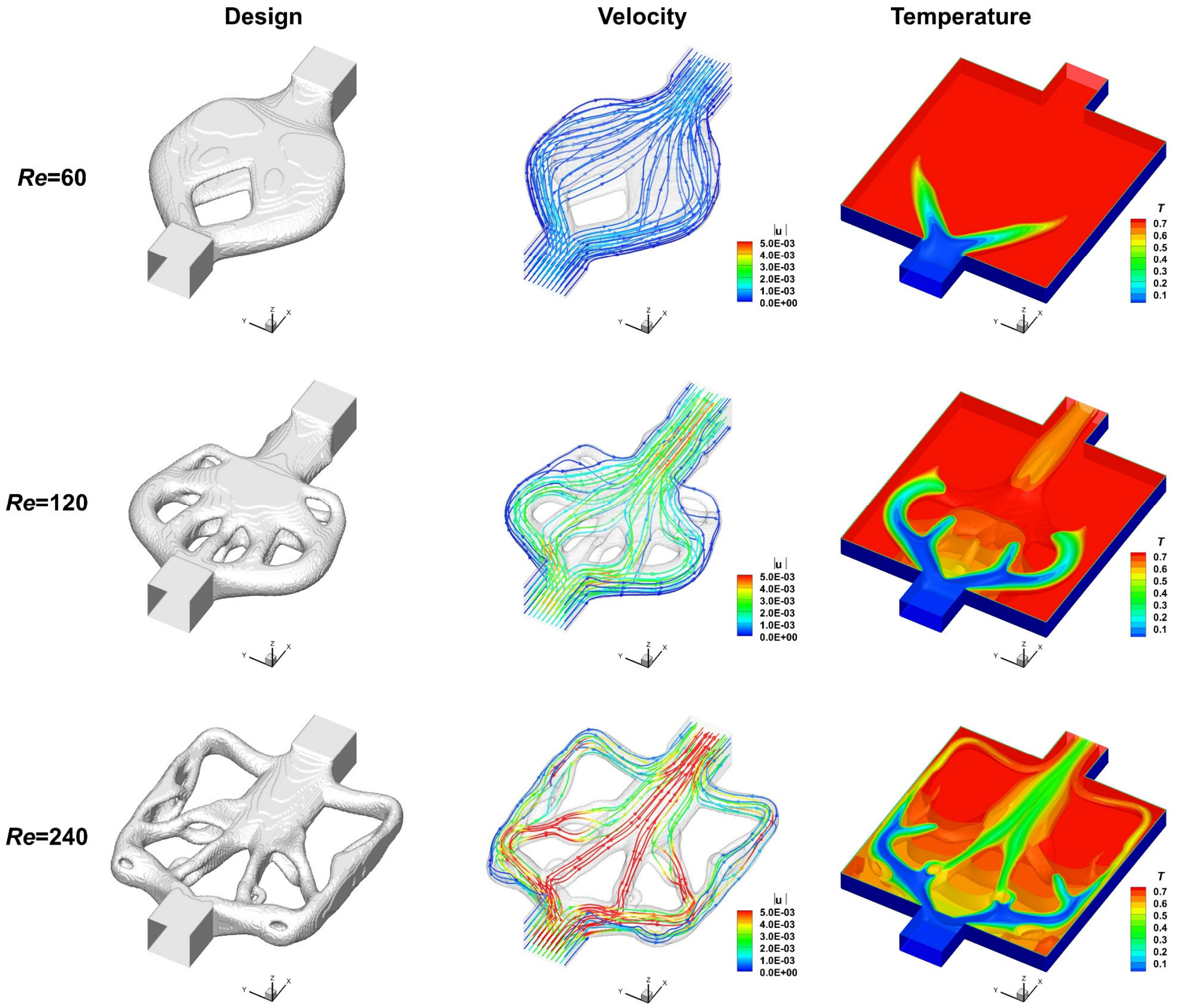

**Figure 19.** Optimized 3D microchannel heat sinks under low fluid volume fraction and the respective velocity and temperature distributions.

3.3.4 Effects of the regularization coefficient

In the base case, the regularization coefficient $\sigma$ is set as a real number, leading to isotropic regularization to the structure during the optimization process. In fact, the regularization coefficient can be arbitrary matrix that enforces anisotropic regularization to the structure, which helps to control the final design more flexibly. Here the regularization coefficient is set as a diagonal matrix $\boldsymbol{\sigma}=\mathrm{diag}(\sigma_{xx}, \sigma_{yy}, \sigma_{zz})$ with $\sigma_{zz}=10\sigma_{xx}$ and $\sigma_{xx}=5\times10^{-3}$. One can expect that the larger regularization coefficient in $z$ direction will bring stronger regularization, which will suppress the structural variation along $z$ direction. Fig. 20 shows the optimized designs under three characteristic $Re$ values. It can be seen that the structural variation along $z$ direction is indeed suppressed, making the optimized designs smoother and more ordered in $z$ direction. In this case, the dimensionless objective values are 0.363, 0.686 and 1.024 respectively, meaning that the stronger regularization in $z$ direction slightly reduces the design freedom and deteriorates the heat transfer performance.

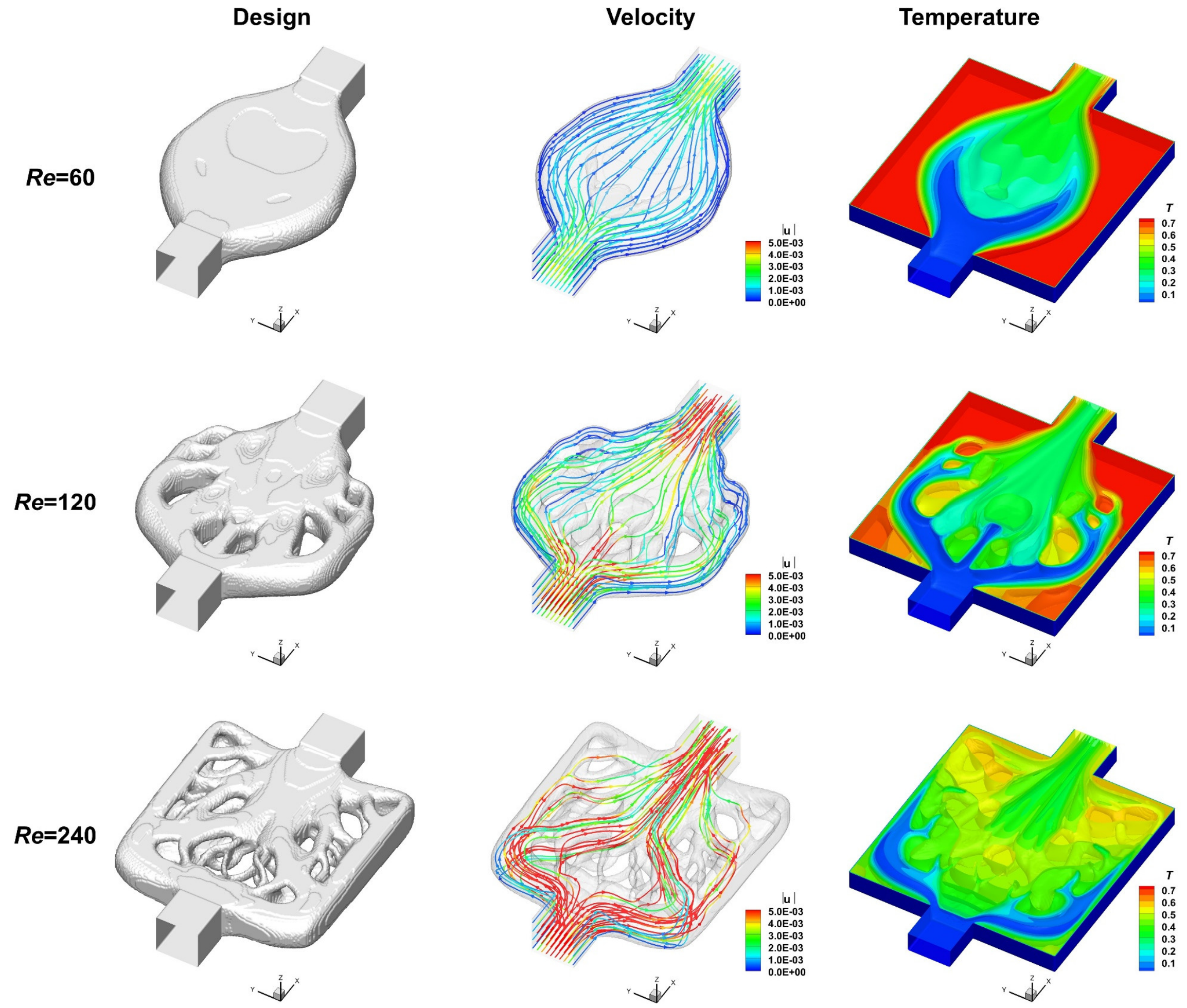

**Figure 20.** Optimized 3D microchannel heat sinks under anisotropic regularization and the respective velocity and temperature distributions.

3.3.5 Comparison to the 2D designs

In the aforementioned cases, the complicated 3D heat sink designs are presented, while in the literature, it is a common situation that the 2D topology optimization is performed and the optimized design is then extruded to generate the 3D design. To reveal the effects of such reduction in dimensionality, the 2D designs are further presented in the present case. Without changing the operating conditions and physical properties, the 2D designs under the same three characteristic *Re* values can be obtained by switching the lattice stencils, and the results are shown in Fig. 21. It can be noted that the optimized designs are also getting more complex, or namely having more individual solid parts, when *Re* increases, and the velocity becomes higher and the temperature in solid gets lower, indicating an enhanced heat transfer. Then these 2D designs are extruded along *z* direction to generate 3D structures, and simulations are performed under the same settings as the base case. The results of 3D simulations are given in Fig. 22. Comparatively, the extruded designs are far less complex than the true 3D optimized designs, since there is no design freedom in *z* direction in the 2D topology optimization. Besides, the effects of the top and bottom walls to the fluid flow and heat transfer are not considered during the 2D topology optimization. Correspondingly, one can expect that the heat transfer performance of the extruded 2D designs will be worse than the 3D optimized designs. In this case, the dimensionless objective values are 0.305, 0.557 and 0.920 respectively, supporting the inference above. However, it is worth mentioning that the better performance of the 3D optimized designs comes at the cost of higher manufacturing expense since the 3D designs are much more complex, and thus compromises should be made when the designers choose the topology optimization methods.

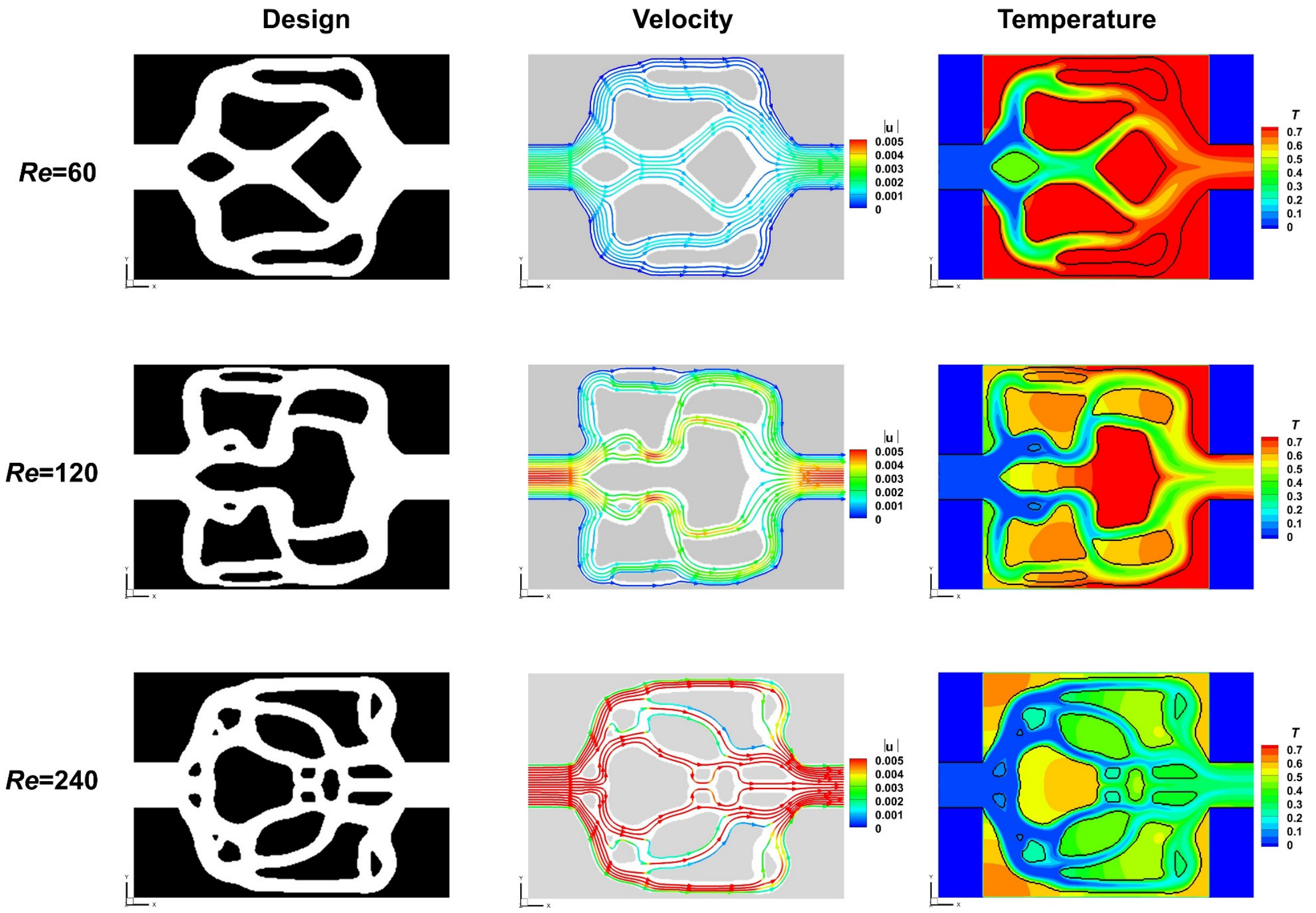

**Figure 21.** Optimized 2D microchannel heat sinks and the respective velocity and temperature distributions.

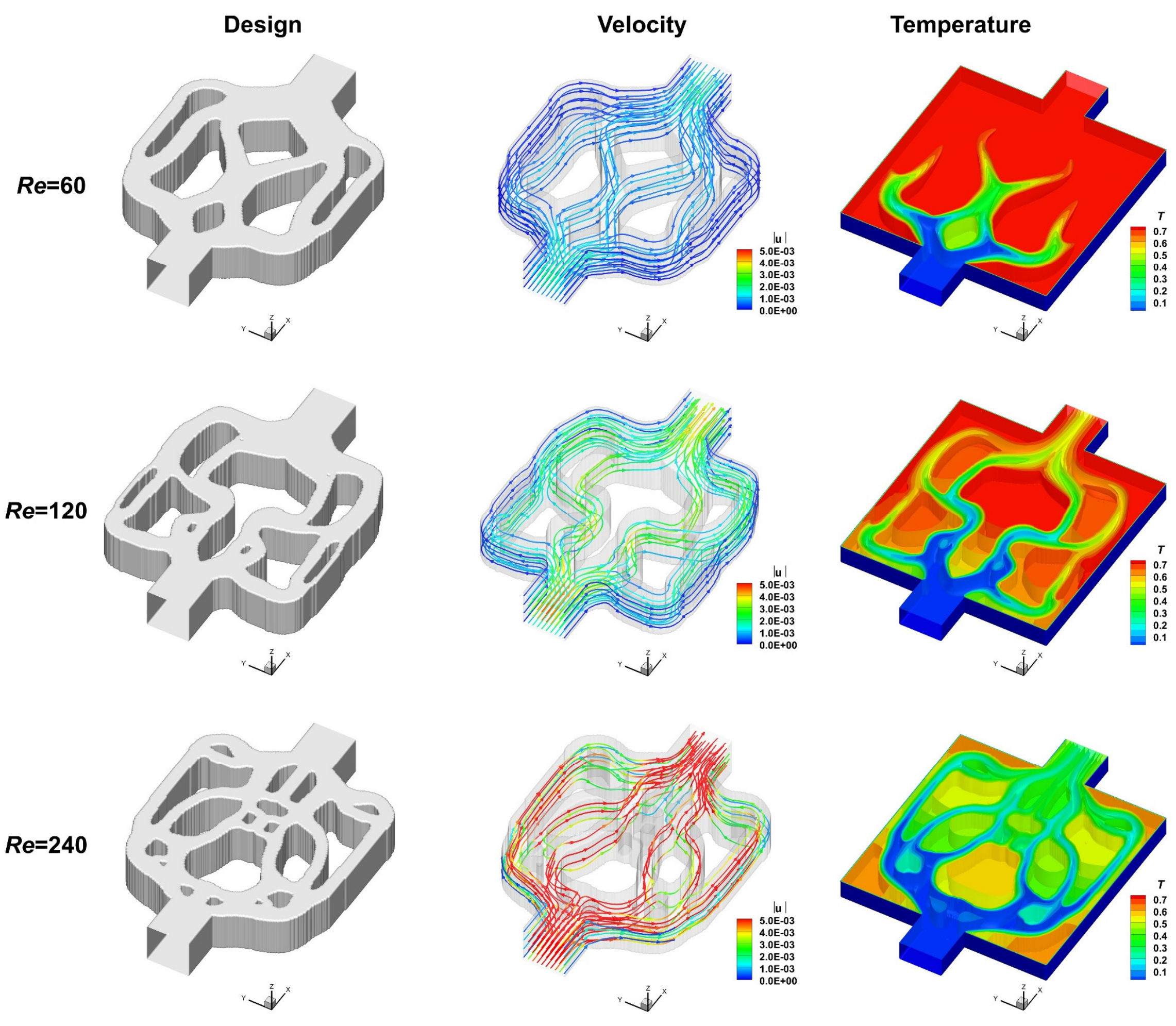

**Figure 22.** Performance verification of the extruded 2D designs.

## 4 Conclusion

This work develops the fully consistent adjoint boundary conditions from the discrete adjoint LBM to improve the theoretical consistency and numerical stability in the continuous adjoint LBM, with both the continuous and discrete adjoint LBMs for open flow systems rigorously derived. It is found that the adjoint collision and adjoint streaming are exactly the same for the two adjoint methods, while the expressions, the objective derivative term and the execution order of the adjoint boundary conditions are totally different. Meanwhile, the continuous adjoint boundary conditions are intrinsically inconsistent or singular, while the discrete adjoint boundary conditions are fully consistent.

The numerical performances of the two derived adjoint models are comprehensively compared by performing sensitivity analysis on the 2D and 3D pipe bend cases. It is found that the continuous adjoint LBM with inconsistent adjoint boundary condition leads to much inferior numerical stability of the adjoint solution and obvious numerical error in sensitivity. The fully consistent discrete adjoint model, however, can always achieve exact sensitivity results and the theoretically highest numerical stability, with a 10 times higher *Re* achieved while without any increase of computational cost. Therefore, the discrete adjoint LBM is a perfect solution to the inconsistency and instability issues in the continuous adjoint LBM. The discrete adjoint LBM is highly recommended for the optimization of open flow systems, while both continuous and discrete adjoint LBM can be acceptable for closed flow systems.

The robust, accurate and efficient adjoint LBM developed above enables the topology optimization of 3D microchannel heat sinks under higher Reynolds number, which are also presented in this work. Esthetic and physically reasonable optimized designs are obtained, demonstrating the necessity and versatility of the presented discrete adjoint LBM. Parametric study finds that reducing the thermal conductivity ratio between solid and fluid leads to a more complex design when Reynolds number is high, and reducing the fluid volume fraction can cause a deteriorated heat transfer due to the weakened convection. Enhancing the regularization in certain direction will slightly reduce the design freedom and deteriorate the heat transfer performance, and the 3D optimized designs can outperform the extruded 2D designs while requiring higher manufacturing cost.

**Acknowledgement**

The authors acknowledge the support of National Nature Science Foundation of China (52376074) and JSPS KAKENHI (Grant No. 23H03799). Ji-Wang Luo also acknowledges the financial support from the China Scholarship Council (CSC).